\documentclass[11pt]{amsart}

\usepackage{amssymb}
\usepackage{amscd}
\usepackage{pb-diagram}
\usepackage{verbatim}
\usepackage{rotating}
\usepackage{graphicx}
\usepackage{enumerate}

\title{Link Concordance and Generalized Doubling Operators}

\author{Tim D. Cochran$^{\dag}$}
\address{Department of Mathematics MS-136, PO Box 1892, Rice University, Houston, TX 77251-1892}
\email{cochran@rice.edu}
\urladdr{http://math.rice.edu/~cochran}

\author{Shelly Harvey$^{\dag\dag}$}
\address{Department of Mathematics MS-136, PO Box 1892, Rice University, Houston, TX 77251-1892}
\email{shelly@rice.edu}
\urladdr{http://math.rice.edu/~shelly}

\author{Constance Leidy}
\address{Department of Mathematics and Computer Science, Wesleyan University, Wesleyan Station, Middletown, CT 06459}
\email{cleidy@wesleyan.edu}
\urladdr{http://cleidy.web.wesleyan.edu}

\thanks{$\dag$The first author was partially supported by NSF grants DMS-0406573 and DMS-0706929.}
\thanks{$\dag\dag$ The second author was partially supported by NSF grant DMS-0539044 and an Alfred P. Sloan Fellowship.}

\newtheorem*{Whitney towers}{Theorem~\ref{Whitney towers}}
\newtheorem*{h-towers}{Theorems ~\ref{half} \& \ref{$(n)$-solvable}}

\newtheorem*{surgery curves}{Theorem~\ref{surgery curves}}
\newtheorem*{cg=0}{Theorem~\ref{vanish}}

\newtheorem{thm}{Theorem}[section]
\newtheorem{lem}[thm]{Lemma}
\newtheorem{cor}[thm]{Corollary}

\newtheorem{prop}[thm]{Proposition}

\theoremstyle{definition}
\newtheorem{defn}[thm]{Definition}
\newtheorem{remark}[thm]{Remark}
\newtheorem{ex}[thm]{Example}

\numberwithin{equation}{section}
\numberwithin{figure}{section}

\newcommand{\Hom}{\operatorname{Hom}}

\newcommand{\G}{\Gamma}

\newcommand{\ra}{\longrightarrow}

\def\x{\times}

\def\ov{\overline}

\def\s{\sigma}
\def\lra{\longrightarrow}

\begin{document}

\begin{abstract} We introduce a technique for showing classical knots and links are not slice. As one application we show that the iterated Bing doubles of many algebraically slice knots are not topologically slice. Some of the proofs do not use the existence of the Cheeger-Gromov bound, a deep analytical tool used by Cochran-Teichner. We define generalized doubling operators, of which Bing doubling is an instance, and prove our nontriviality results in this more general context. Our main examples are boundary links that cannot be detected in the algebraic boundary link concordance group.
\end{abstract}

\maketitle

%%%%%%%%%%%%%%%%%%%%   Start of main body of article

\section{Introduction}\label{sec:Introduction}

A \textbf{link} $L=\{K_1,...,K_m\}$ of $m$-components is an ordered  collection of $m$ oriented circles disjointly embedded in $S^3$. A \textbf{knot} is a link of one component. A  \textbf{topological slice link} (abbreviated as \textbf{slice} in this paper) is a link whose components bound a disjoint union of $m$ $2$-disks topologically and locally flatly embedded in $B^4$. The question of which links are slice links lies at the heart of the topological classification of $4$-dimensional manifolds.

The connected sum operation gives the set of all knots, modulo slice knots, the structure of an abelian group, called the \emph{topological knot concordance group $\mathcal{C}$}, which is a quotient of its smooth analogue. For links one must consider \emph{string} links to get a well-defined group structure, and this operation is not commutative ~\cite{Le7}. This group is called the $m$-component \emph{string link concordance group}. We applied our techniques to knot concordance in ~\cite{CHL3}\cite{CHL1}. This paper gives new information about link concordance. All of the results here (except Example~\ref{ex:nonzerofirst-ordersigs}) were first announced in ~\cite{CHL1} and appeared in the unpublished preprint ~\cite{CHL1A} under a different title (that preprint was later split into two papers, the present paper being one). We employ the Cheeger-Gromov von Neumann $\rho$-invariants and higher-order Alexander modules that were introduced in ~\cite{COT}. Our new technique is to expand upon previous results of Leidy concerning higher-order Blanchfield forms for arbitrary $3$-manifolds ~\cite{Lei3}~\cite{Lei1}. This is used to show that certain elements of $\pi_1$ of a  link exterior cannot lie in the kernel of the map into $\pi_1$ of a slice disk exterior. We also employ results of Harvey on the \emph{torsion-free derived series of groups} ~\cite{Ha2}, and recent results of Cochran-Harvey on versions of Dwyer's Theorem for the derived series ~\cite{CH2}. We note that the construction of examples is in the smooth category so that we actually also prove the corresponding statements about smooth link concordance.

\begin{figure}[htbp]
\begin{center}
\setlength{\unitlength}{1pt}
\begin{picture}(112,109)
\put(0,0){\includegraphics{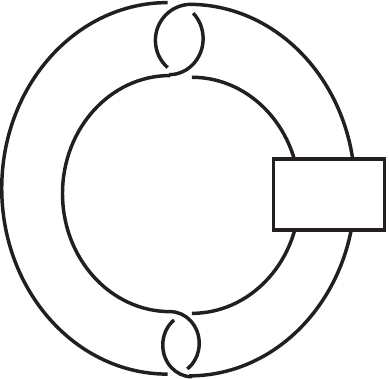}}
\put(88,50){$K$}
\put(120,54){=~~$BD(K)$}
\end{picture}
\end{center}
\caption{Bing double of $K$}\label{fig:Bingdouble}
\end{figure}

Natural families of links have been considered. In particular, if $K$ is any knot then the \textbf{Bing double} of $K$, $BD(K)$ is the $2$-component link shown in Figure~\ref{fig:Bingdouble}.

If $K$ is slice then it is easy to see that $BD(K)$ is a slice link. A natural question is whether or not the converse is true. It was shown by Harvey that if the Bing double (or even an iterated Bing double) of $K$ is topologically slice then the integral over the circle of the Levine signatures of $K$ is zero ~\cite[Corollary 5.6]{Ha2}. It was shown by Cimasoni that if $BD(K)$ is a \emph{boundary slice} link then $K$ is algebraically slice ~\cite{Ci}. Subsequently (and after ~\cite{CHL1}) it was shown by Cha-Livingston-Ruberman that if $BD(K)$ is a slice link then $K$ must be an algebraically slice knot ~\cite{CRL}. Here we address the questions: If $K$ is algebraically slice then does it follow that $BD(K)$ is a topological slice link? What about for iterated Bing doubles? We answer these questions in the negative by showing that certain higher-order signatures of $K$ offer further obstructions. For example, in Section~\ref{sec:Bingdoubles} we define \emph{first-order signatures} of $K$, akin to Casson-Gordon invariants, and show that the first-order signatures of $K$, like the ordinary signatures, obstruct any \emph{iterated Bing double} of $K$ from being a slice link. This improves on Harvey's theorem.

\newtheorem*{thm:Bingdouble}{Theorem~\ref{thm:Bingdouble}}
\begin{thm:Bingdouble} Let $K$ be an arbitrary knot. If some iterated Bing double of $K$ is topologically slice in a rational homology $4$-ball then one of the first-order signatures of $K$ is zero.
\end{thm:Bingdouble}

\begin{figure}[htbp]
\begin{center}
\setlength{\unitlength}{1pt}
\begin{picture}(143,151)
\put(0,0){\includegraphics{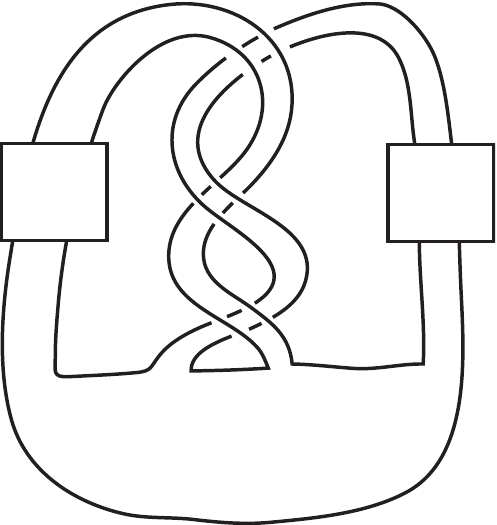}}
\put(120,94){$J_0$}
\put(10,94){$J_0$}
\put(-30,94){$J_1\equiv$}
\end{picture}
\end{center}
\caption{Algebraically Slice Knots $J_1$}\label{fig:algsliceexs}
\end{figure}
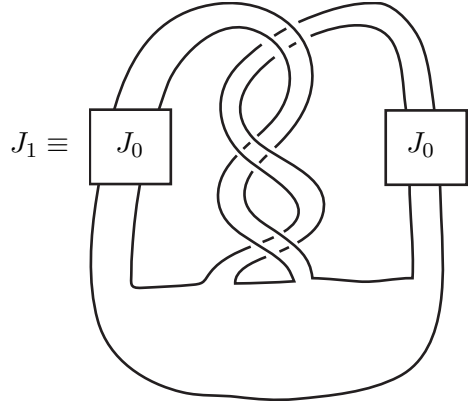

\begin{figure}[htbp]
\begin{center}
\begin{picture}(138,135)
\put(0,0){\includegraphics{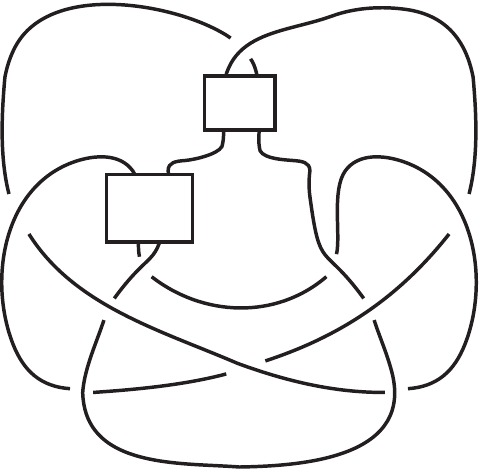}}
\put(-30,70){$K_1~=$}
\put(61,103){$K_0$}
\put(37,73){$K_0$}
\end{picture}
\end{center}
\caption{Algebraically slice knots $K_1$}\label{fig:first8_9}
\end{figure}

\begin{figure}[htbp]
\setlength{\unitlength}{1pt}
\begin{center}
\begin{picture}(159,101)
%\put(0,0){\includegraphics{figure8_150pt.pdf}}
\put(0,0){\includegraphics{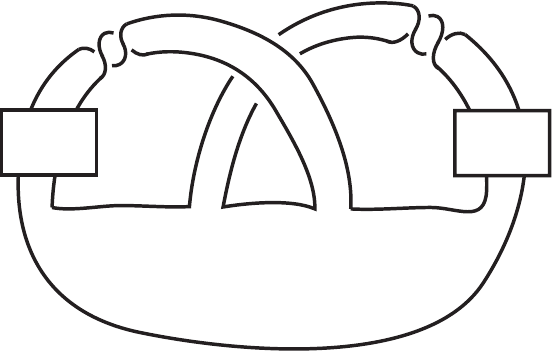}}
\put(-25,50){$E=$}
\put(8,57){$E_0$}
\put(135,57){$E_0$}
%\put(5,50){$E_1=$}
%\put(46,57){$E_0$}
%\put(178,57){$E_0$}
\end{picture}
\end{center}
\caption{}\label{fig:examplehighersigsfigeight2}
\end{figure}

For example, for the algebraically slice knots $J_1$ of Figure~\ref{fig:algsliceexs} and $K_1$ of Figure~\ref{fig:first8_9}, the first order signatures are related to classical signatures of $J_0$ and $K_0$ respectively, and similarly for the knots $E_1$ as in Figure~\ref{fig:examplehighersigsfigeight2}, which are of order 2 in the algebraic concordance group.

\begin{def}\label{def:rho0} For a knot $K$ in $S^3$ let $\boldsymbol{\rho_0(K)}$ denote the integral over the circle of the classical Levine signature function of $K$ (normalize so that the length of the circle is $1$).
\end{def}

On the examples above, Theorem~\ref{thm:Bingdouble} takes the following nice form:

\newtheorem*{cor:Bingdouble2}{Corollary~\ref{cor:Bingdouble2}}
\begin{cor:Bingdouble2}If $K_1$ is an algebraically slice knot of Figure~\ref{fig:first8_9} and if some iterated Bing double of $K_1$ is slice in a rational homology ball then $\rho_0(K_0)=0$. If $E_1$ is a knot as in Figure~\ref{fig:examplehighersigsfigeight2} and some iterated Bing double of $E_1$ is slice in a rational homology ball then $\rho_0(E_0)=0$.
\end{cor:Bingdouble2}

\newtheorem*{cor:Bingdouble}{Corollary~\ref{cor:Bingdouble}}
\begin{cor:Bingdouble} There is a constant $C$ such that, if  $\rho_0(J_0)\notin \{0,C\}$ then no iterated Bing double of the algebraically slice knot in Figure~\ref{fig:algsliceexs}  is slice in a rational homology ball.
\end{cor:Bingdouble}

It is well known that the Bing double of any knot is a boundary link and that the Bing double of any algebraically slice knot is zero in the algebraic boundary link concordance group (proofs can be found in ~\cite[Proposition 1.1, Theorem 2.1 (i)]{Ci}). Thus the examples above are boundary links that cannot be detected in the algebraic boundary link concordance group.

We remark that recent work of Cha shows that even many amphichiral knots have non-slice Bing doubles ~\cite{Cha4}. Amphichiral knots cannot be handled by the techniques in the present paper.

\begin{figure}[htbp]
\setlength{\unitlength}{1pt}
\begin{center}
\begin{picture}(143,151)
\put(0,0){\includegraphics{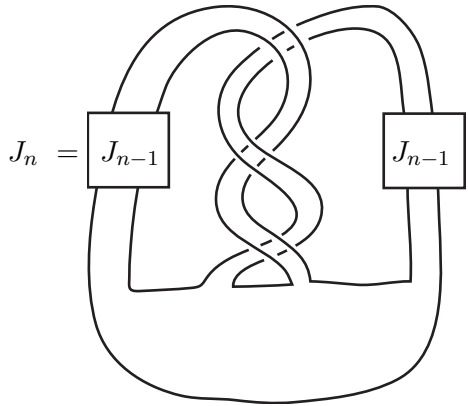}}
\put(05,92){$J_{n-1}$}
\put(-30,92){$J_{n}~=$}
\put(115,92){$J_{n-1}$}
%\put(30,10){\includegraphics{family_150pt.pdf}}
%\put(151,101){$J_0$}
%\put(40,100){$J_0$}
%\put(0,98){$J_1\equiv$}
\end{picture}
\end{center}
\caption{The recursive family $J_{n+1}, n\geq 0$}\label{fig:family}
\end{figure}

We have similar results for iterated Bing doubles of even more subtle knots. For example, consider the (recursively-defined) family $J_n, ~n>0$, of Figure~\ref{fig:family} whose members are not only algebraically slice but also have vanishing Casson-Gordon invariants for every $n>1$. An $n^{th}$-order higher-order signature of $J_n$ obstructs the iterated Bing doubles of $J_n$ from being slice links. Moreover these iterated Bing doubles give non-trivial examples of links that lie deeper and deeper in the Cochran-Orr-Teichner filtration of the set of concordance classes of links
$$
\cdots \subseteq \mathcal{F}_{n} \subseteq \cdots \subseteq
\mathcal{F}_1\subseteq \mathcal{F}_{0.5} \subseteq \mathcal{F}_{0} \subseteq \mathcal{C}.
$$
This filtration, first defined in ~\cite[Sections 7,8]{COT}, is reviewed in Section~\ref{sec:appendix}. Recall that a link in $\mathcal{F}_{n}$ is called \textbf{$(n)$-solvable}. The knot $J_n$ is $(n)$-solvable but not necessarily $(n+1)$-solvable.

\newtheorem*{cor:BingdoubleJnK}{Corollary~\ref{cor:BingdoubleJnK}}
\begin{cor:BingdoubleJnK} For any $n$ there is a constant $C$ such that for any knot $J_0$ with Arf invariant zero and $|\rho_0(J_0)|>C$ the Bing double of $J_{n-1}$ is $(n)$-solvable but not slice nor even $(n+1)$-solvable.
\end{cor:BingdoubleJnK}

\newtheorem*{cor:iteratedbing}{Corollary~\ref{cor:iteratedbing}}
\begin{cor:iteratedbing} Suppose $k$ and $n$ are positive integers. Then there is a constant $C$ such that for any knot $J_0$ with Arf($J_0)=0$ and $|\rho_0(J_0)|>C$, the $k$-fold iterated Bing double of $J_{n-k}$ is $(n)$-solvable but not slice nor even $(n+1)$-solvable.
\end{cor:iteratedbing}

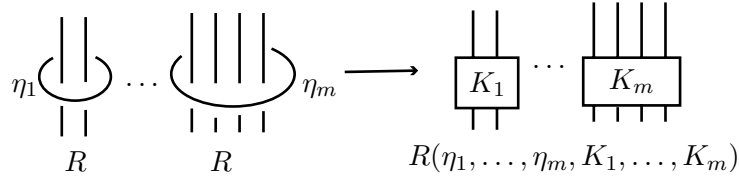
\begin{figure}[htbp]
\setlength{\unitlength}{1pt}
\begin{center}
\begin{picture}(243,61)
\put(-10,27){$\eta_1$} \put(100,27){$\eta_m$} \put(33,29){$\dots$}
\put(187,36){$\dots$} \put(163,27){$K_1$} \put(216,28){$K_m$}
\put(140,-1){$R(\eta_1,\dots,\eta_m,K_1,\dots,K_m)$}
\put(10,-3){$R$} \put(65,-3){$R$}
\put(0,10){\includegraphics{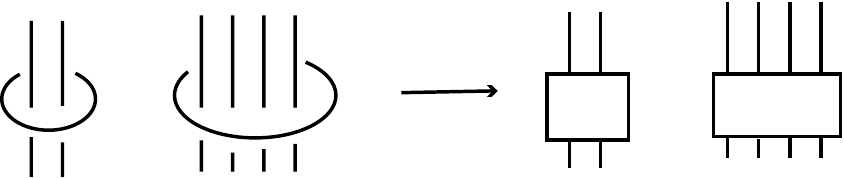}}
\end{picture}
\end{center}
\caption{$R(\eta_1,\dots,\eta_m,K_1,\dots,K_m)$:
Infection of $R$ by $K_i$ along $\eta_i$}\label{fig:infection}
\end{figure}

The specific families of links of Figure~\ref{fig:Bingdouble} are important because of their simplicity. However, they are merely particular instances of a more general doubling phenomenon to which our techniques may be applied. In order to state these results, we review a method we will use to construct examples. Let $R$ be a link in $S^3$ and let $\{\eta_1,\eta_2,\ldots,\eta_m\}$ be an oriented trivial link in $S^3$ that misses $R$ and  bounds a collection of disks that meet $R$ transversely. Suppose $\{K_1,K_2,\ldots,K_m\}$ is an $m$-tuple of auxiliary knots. Let $R(\eta_1,\ldots,\eta_m,K_1,\ldots,K_m)$ denote the result of the operation pictured in Figure~\ref{fig:infection}, that is, for each $\eta_i$, take the embedded disk in $S^3$ bounded by $\eta_i$; cut off $R$ along the disk; grab the cut strands, tie them into the knot $K_i$ (with no twisting) and reglue as shown in Figure~\ref{fig:infection}.

\begin{figure}[htbp]
\setlength{\unitlength}{1pt}
\begin{center}
\begin{picture}(113,109)
\put(0,0){\includegraphics{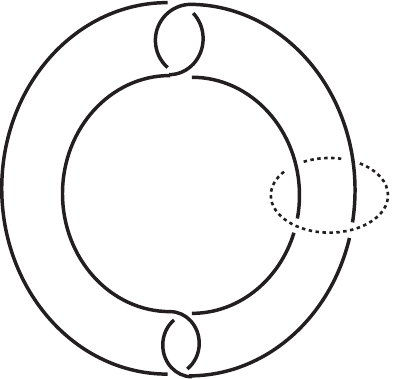}}
\put(120,50){$\alpha$}
\put(137,54){}
\end{picture}
\end{center}
\caption{Bing double of $K$ is infection on the trivial link along $\alpha$} using $K$ \label{fig:bingeta}
\end{figure}

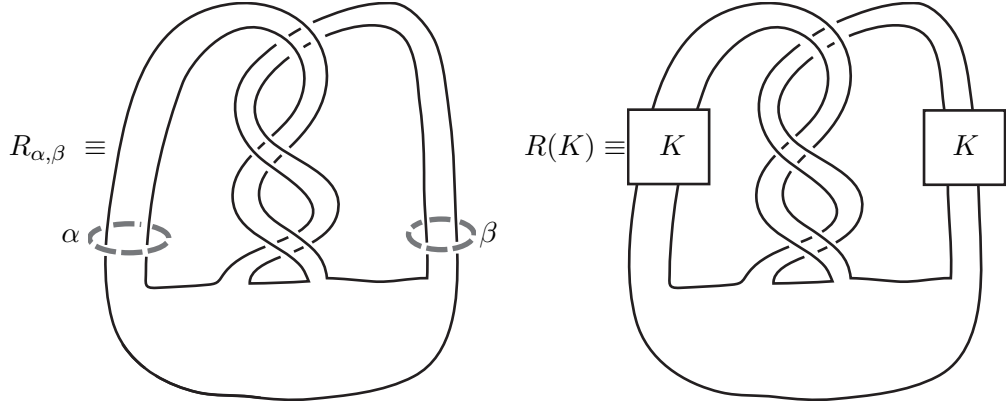
\begin{figure}[htbp]
\setlength{\unitlength}{1pt}
\begin{center}
\begin{picture}(357,151)
\put(20,0){\includegraphics{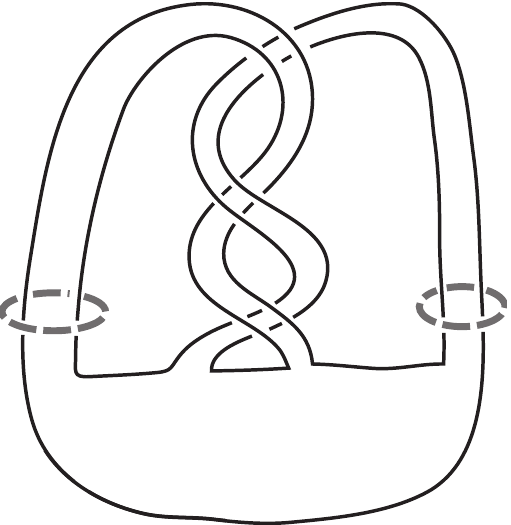}}
\put(224,0){\includegraphics{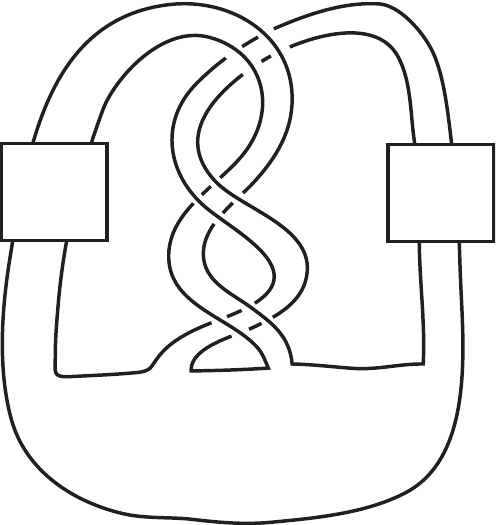}}
\put(10,60){$\alpha$}
\put(168,60){$\beta$}
\put(236,93){$K$}
\put(347,93){$K$}
\put(185,93){$R(K)\equiv$}
\put(-10,93){$R_{\alpha,\beta}~\equiv$}
\end{picture}
\end{center}
\caption{$R$-doubling}\label{fig:Rdoubling}
\end{figure}

\noindent We will call this the result of \textbf{infection performed on the link $\boldsymbol{R}$ using the infection knots $\boldsymbol{K_i}$ along the curves $\boldsymbol{\eta_i}$}. This construction can also be
described in the following way. For each $\eta_i$, remove a tubular neighborhood of $\eta_i$ in $S^3$ and glue in the exterior of a tubular neighborhood of $K_i$ along their common boundary, which is a
torus, in such a way that the longitude of $\eta_i$ is identified with the meridian of $K_i$ and the meridian of $\eta_i$ with the reverse of the longitude of $K_i$. The resulting space can be seen to be homeomorphic to $S^3$ and the image of $R$ is the new link. In the case that $m=1$ this is the same as the classical satellite construction. In general it can be considered to be a `generalized satellite construction', widely utilized in the study of knot concordance. In the case that $m=1$ and $lk(\eta,R)=0$ it is precisely the same as forming a satellite of $J$ with winding number zero. This yields an operator
$$
R_{\eta}:\mathcal{C}\to \mathcal{C}^k.
$$
where $\mathcal{C}^k$ is the set of concordance classes of $k$-component links. For general $m$ with $lk(\eta_i,R)=0$, it should be considered as a \textbf{generalized doubling operator}, $R_{\eta_i}$, parameterized by $(R,\{\eta_i\})$. If, for simplicity, we assume that all input knots are identical then such an operator is a function
$$
R_{\eta_i}:~\mathcal{C}\to \mathcal{C}^k.
$$
Bing-doubling is an example of this ($m=1$) as suggested by Figure~\ref{fig:bingeta}.
Another primary example is the ``$9_{46}$-doubling'' operation of going from the left-hand side of Figure~\ref{fig:Rdoubling} to the right-hand side.
Here $R$ is the $9_{46}$ knot and $\{\eta_1,\eta_2\}=\{\alpha,\beta\}$ are as shown on the left-hand side of Figure~\ref{fig:Rdoubling}. The image of a knot $K$ under the operator $R_{\alpha,\beta}$ is denoted by $R(K)$ and is shown on the right-hand side of Figure~\ref{fig:Rdoubling}. Note that our previously defined knot $J_1$ is the same as $R(J_0)$.

\noindent Most of the results of this paper concern to what extent these functions are injective. The point is that, because of the condition on ``winding numbers,'' $lk(\eta_i,R)=0$,
if $R$ is a slice link, the images of such operators $R$ contain only links for which the classical Seifert-matrix-type invariants vanish. Moreover these operators respect the COT filtration.

\newtheorem*{lem:nsolv}{Lemma~\ref{lem:nsolv}}
\begin{lem:nsolv} If $R$ is a slice link and $\eta_i\in \pi_1(S^3-R)^{(p)}$ then the operator $R_{\eta_i}$ satisfies
$$
R_{\eta_i}(\mathcal{F}_{q})\subset \mathcal{F}_{p+q}.
$$
\end{lem:nsolv}

Thus iterations of these operators, \textbf{iterated generalized doubling}, produce increasingly subtle links. More generally let us define an \textbf{$n$-times iterated generalized doubling} to be precisely such a composition of operators using possibly different slice links $R_j$, and different curves $\eta_{j1},\dots,\eta_{jm_j}$. For example the knot $J_n$ of Figure~\ref{fig:family} is obtained from $J_0$ by applying $R\circ\dots\circ R$ where $R=R_{\alpha,\beta}$ is as in Figure~\ref{fig:Rdoubling}. Then, generalizing Corollary~\ref{cor:iteratedbing}, our method establishes:

\newtheorem*{thm:mainlink3}{Theorem~\ref{thm:mainlink3}}
\begin{thm:mainlink3} Suppose $T$ is a slice link, $\alpha$ is an unknotted circle in $S^3-T$ that represents an element in $\pi_1(S^3-T)^{(k)}$ but not in $\pi_1(M_T)^{(k+1)}_H$. Suppose for each $j$, $1\leq j\leq n-k$, $R_j$ is a slice knot, $\{\eta_{j1},\dots,\eta_{jm_j}\}$ is a trivial link of circles in $S^3-R_j$ with the property that the submodule of the classical Alexander polynomial of $R_j$ generated by $\{\eta_{j1},\dots,\eta_{jm_j}\}$ contains elements $x,y$ such that $\mathcal{B}\ell_0^j(x,y)\neq 0$, where $\mathcal{B}\ell_0^j$ is the Blanchfield form of $R_j$. Finally suppose that Arf($K$)$=0$. Then the result, $L(K)\equiv T_{\alpha}\circ R_{n-k}\circ\dots\circ R_1(K)$, of the iterated generalized doubling (applied to $K$) lies in $\mathcal{F}_{n}$ and there is a constant $C$ (independent of $K$), such that if $|\rho_0(K)|>C$, then $L(K)$ is of infinite order in the topological concordance group (moreover no multiple lies in $\mathcal{F}_{n+1}$).
\end{thm:mainlink3}

\section{Higher-Order Signatures and How to Calculate Them}\label{signatures}

In this section we review the von Neumann $\rho$-invariants and explain to what extent they are concordance invariants. We also show how to calculate them for knots or links that are obtained from the infections defined in Section~\ref{sec:Introduction}.

The use of variations of Hirzebruch-Atiyah-Singer signature defects associated to covering spaces is a theme common to most of the work in the field of knot and link concordance since the 1970's. In particular, Casson and Gordon initiated their use in cyclic covers ~\cite{CG1}~\cite{CG2}; Farber, Levine and Letsche initiated the use of signature defects associated to general (finite) unitary representations ~\cite{Let}~\cite{L6}; and Cochran-Orr-Teichner initiated the use of signatures associated to the left regular representations ~\cite{COT}. See ~\cite{Fr2} for a beautiful comparison of these approaches in the metabelian case.

Given a compact, oriented 3-manifold $M$, a discrete group $\G$, and a representation $\phi : \pi_1(M)
\to \G$, the \emph{von Neumann
$\boldsymbol{\rho}$-invariant} was defined by Cheeger and Gromov by choosing a Riemannian metric and using $\eta$-invariants associated to $M$ and its covering space induced by $\phi$. It can be thought of as an oriented homeomorphism invariant associated to an arbitrary regular covering space of $M$ ~\cite{ChGr1}. If $(M,\phi) = \partial
(W,\psi)$ for some compact, oriented 4-manifold $W$ and $\psi : \pi_1(W) \to \G$, then it is known that $\rho(M,\phi) =
\s^{(2)}_\G(W,\psi) - \s(W)$ where $\s^{(2)}_\G(W,\psi)$ is the
$L^{(2)}$-signature (von Neumann signature) of the intersection form defined on
$H_2(W;\mathbb{Z}\G)$ twisted by $\psi$ and $\sigma(W)$ is the ordinary
signature of $W$ ~\cite{LS}. In the case that $\G$ is a poly-(torsion-free-abelian) group (abbreviated \textbf{PTFA group} throughout), it follows that $\mathbb{Z}\G$ is a right Ore domain that embeds into its (skew) quotient field of fractions $\mathcal{K}\G$ ~\cite[pp.591-592, ~Lemma 3.6ii p.611]{P}. In this case $\s^{(2)}_\G$  is a function of the Witt class of the equivariant intersection form on $H_2(W;\mathcal{K}\G)$ ~\cite[Section 5]{COT}. In the special case that this form is non-singular (such as $\beta_1(M)=1$), it can be thought of as a homomorphism from
$L^0(\mathcal{K}\G)$ to $\mathbb{R}$.

All of the coefficient systems $\G$ in this paper will be of the form $\pi/\pi^{(n)}_r$ where $\pi$ is the fundamental group of a space (usually a $4$-manifold) and $\pi^{(n)}_r$ is the $n^{th}$-term of the \textbf{rational derived series}. The latter was first considered systematically by Harvey. It is defined by
$$
\pi^{(0)}_r\equiv \pi,~~~ \pi^{(n+1)}_r\equiv \{x\in \pi^{(n)}_r|\exists k\neq 0, x^k\in [\pi_r^{(n)},\pi_r^{(n)}]\}.
$$
Note that $n^{th}$-term of the usual derived series $\pi^{(n)}$ is contained in the $n^{th}$-term of the rational derived series. For free groups and knot groups, they coincide. It was shown in ~\cite[Section 3]{Ha1} that $\pi/\pi^{(n)}_r$ is a PTFA group.

The utility of the von Neumann signatures lies in the fact that they obstruct knots from being slice knots. It was shown in ~\cite[Theorem 4.2]{COT} that, in certain situations, higher-order von Neumann signatures vanish for slice knots, generalizing the classical result of Murasugi and the results of Casson-Gordon. That proof fails for links, but the extension was later accomplished by Harvey (there is an extra obstruction). Moreover, Cochran-Orr-Teichner defined a filtration on knots and links and showed that certain higher-order signatures obstructed a knot's lying in a certain term of the filtration. Harvey also extended this to links. In this section we state the needed results for \emph{slice} knots and links. For those readers interested primarily in what links are slice, this suffices. For those readers interested in the $(n)$-solvable filtration, we have included in Section~\ref{sec:appendix} a review of this filtration as well as important results about vanishing of $\rho$ invariants (some new). Such readers would be advised to read Section~\ref{sec:appendix} after finishing the present section.

First, we recall the theorem of Cochran-Orr-Teichner for knots.

\begin{thm}[Cochran-Orr-Teichner~{\cite[Theorem 4.2]{COT}}]
\label{thm:oldsliceobstr} If a knot $K$ is topologically slice in a rational homology $4$-ball and
$\phi:\pi_1(M_K)\to \G$ is a PTFA coefficient system that extends to the fundamental group of the exterior of the slicing disk, then $\rho(M_K,\phi)=0$.
\end{thm}

The analogous result for links has only recently appeared, although it is implicit in and follows from the results of ~\cite{Ha2}.

\begin{thm}[Corollary of Cochran-Harvey~{\cite[Theorem 4.9, Proposition 4.11]{CH2}}]
\label{thm:linksliceobstr}If a link $L$ is topologically slice in a rational homology $4$-ball and
$\phi:\pi_1(M_L)\to \G$ is a PTFA coefficient system that extends to the fundamental group of the exterior of the slicing disks, then $\rho(M_L,\phi)=0$.
\end{thm}

Some other useful properties of von Neumann
$\rho$-invariants are given below. One can find
detailed explanations of most of these in \cite[Section 5]{COT}. The last property, that for a fixed $3$-manifold, the set $\{\rho(M,\phi)\}$ is bounded above and below, is an analytical result of Cheeger and Gromov that we use in some (but not all) of our results here.

\begin{prop}
\label{prop:rho invariants}Let $M$ be a closed, oriented $3$-manifold and $\phi : \pi_1(M) \to \G$ as above.
\begin{itemize}
\item [(1)] If $(M,\phi) = \partial (W,\psi)$ for some compact
oriented 4-manifold $W$ such that the equivariant intersection form on $H_2(W;\mathcal{K}\G)/j_*(H_2(\partial W;\mathcal{K}\G))$ admits a half-rank summand on
which the form vanishes, then
$\s^{(2)}_\G(W,\psi)=0$ (see ~\cite[Lemma 3.1 and Remark 3.2]{Ha2} for a proper explanation of this for manifolds with $\beta_1>1$). Thus if  $\s(W) = 0$ then $\rho(M,\phi) = 0.$

\item [(2)] If $\phi$ factors through $\phi' : \pi_1(M) \to \G'$ where
$\G'$ is a subgroup of $\G$, then $\rho(M,\phi') = \rho(M,\phi)$.
%\mathbf{}
\item [(3)] If $\phi$ is trivial (the zero map), then $\rho(M,\phi) = 0$.

\item [(4)] If $M=M_K$ is zero surgery on a knot $K$ and $\phi:\pi_1(M)\to \mathbb{Z}$ is the abelianization, then $\rho(M,\phi)$ is equal to $\rho_0(K)$ ~\cite[Prop. 5.1]{COT2}.

\item [(5)] (Cheeger-Gromov ~\cite{ChGr1}) Given $M$, there is a positive constant $\boldsymbol{C_M}$, the \textbf{Cheeger-Gromov constant} of $M$, such that for every $\phi$
$$
|\rho(M,\phi)|<C_M.
$$
\end{itemize}
\end{prop}

The following elementary lemma reveals the additivity of the $\rho$-invariant under infection. It is only slightly more general than \cite[Proposition 3.2]{COT2}. The use of a Mayer-Vietoris sequence to analyze the effect of a satellite construction on signature defects is common to essentially all of the previous work in this field (see for example ~\cite{Lith1}).

Suppose $L=R(\eta_i,K_i)$ is obtained by infection as described in Section~\ref{sec:Introduction}. Let the zero surgeries on $R$, $L$, and $K_i$ be denoted $M_R$ $M_L$, $M_i$ respectively. Suppose $\phi:\pi_1(M_L)\to \G$ is a map to an arbitrary PTFA group $\G$ such that, for each $i$, $\ell_i$, the longitude of $K_i$, lies in the kernel of $\phi$. Since $S^3-K_i$ is a submanifold of $M_L$, $\phi$ induces  a map on $\pi_1(S^3-K_i)$. Since $l_i$ lies in the kernel of $\phi$, this map extends uniquely to a map that we call $\phi_i$  on $\pi_1(M_i)$. Similarly, $\phi$ induces a map on $\pi_1(M_R-\coprod \eta_i)$. Since $M_R$ is obtained from $(M_R-\coprod \eta_i)$ by adding $m$ $2$-cells along the meridians of the $\eta_i$, denoted $\mu(\eta_i)$, and $m$ $3-$cells, and since $\mu(\eta_i)=l_i^{-1}$ and $\phi_i(l_i)=1$, $\phi$ extends uniquely to $\phi_R$. Thus $\phi$ induces unique maps $\phi_i$ and $\phi_R$ on $\pi_1(M_i)$ and $\pi_1(M_R)$ (characterized by the fact that they agree with $\phi$ on $\pi_1(S^3-K_i)$ and $\pi_1(M_R-\coprod \eta_i)$ respectively).

There is a very important case when the hypothesis above that $\phi(\ell_i)=1$ is always satisfied. Namely suppose $\G^{(n+1)}=1$ and $\eta_i\in \pi_1(M_R)^{(n)}$.  Since a longitudinal push-off of $\eta_i$, called $\ell_{\eta_i}$ or $\eta_i^+$, is isotopic to $\eta_i$ in the solid torus $\eta_i\times D^2\subset M_R$, $\ell_{\eta_i}\in \pi_1(M_R)^{(n)}$ as well. By ~\cite[Theorem 8.1]{C} or ~\cite{Lei3} it follows that $\ell_{\eta_i}\in \pi_1(M_L)^{(n)}$. Since $\mu_i$, the meridian of $K_i$, is identified to $\ell_{\eta_i}$, $\mu_i \in \pi_1(M_L)^{(n)}$ so $\phi(\mu_i)\in \G^{(n))}$  for each $i$.  Thus $\phi_i(\pi_1(S^3- K_i)^{(1)})\subset\G^{(n+1)}=\{e\}$ and in particular the longitude of each $K_i$ lies in the kernel of $\phi$.

\begin{lem}[{\cite[Lemma 2.3]{CHL3}}]\label{lem:additivity} In the notation of the two previous paragraphs (assuming $\phi(\ell_i)=0$ for all $i$),
$$
\rho(M_L,\phi) - \rho(M_R,\phi_R) = \sum^m_{i=1}\rho(M_i,\phi_i).
$$
Moreover if $\pi_1(S^3-K_i)^{(1)}\subset$ kernel($\phi_i$) then either $\rho(M_i,\phi_i)=\rho_0(K_i)$, or $\rho(M_i,\phi_i)=0$, according as $\phi_R(\eta_i)\neq 1$ or $\phi_R(\eta_i)= 1$. Specifically, if $\G^{(n+1)}=1$ and $\eta_i\in \pi_1(M_R)^{(n)}$ then this is the case.
\end{lem}

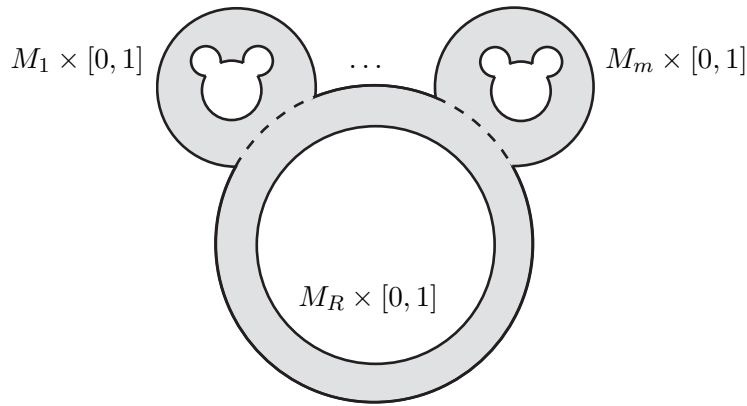
\begin{figure}[htbp]
\setlength{\unitlength}{1pt}
\begin{center}
\begin{picture}(150,150)
\put(0,0){\includegraphics{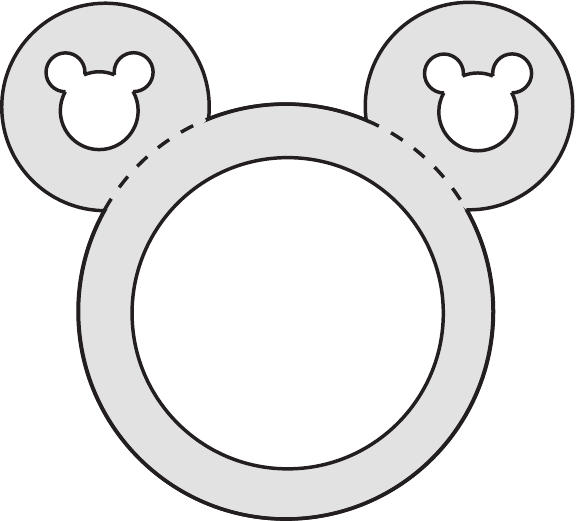}}
\put(54,37){$M_R\times [0,1]$}
\put(-55,127){$M_1\times [0,1]$}
\put(170,127){$M_m\times [0,1]$}
\put(73,127){$\dots$}
\end{picture}
\end{center}
\caption{The cobordism $E$}\label{fig:mickey}
\end{figure}

We will now sketch the proof since we need, independently, several elements of that proof. Most importantly, there is a cobordism that relates the zero surgeries of the original link, the link achieved through infection(s) and the zero surgeries on the infecting knots. Let $E$ be the $4$-manifold obtained from $M_R\times [0,1] \coprod -M_i\times [0,1]$ by identifying, for each $i$, the copy of $\eta_i\times D^2$ in $M_R\times \{1\}$ with the tubular neighborhood of $K_i$ in $M_i\times \{0\}$ as in Figure~\ref{fig:mickey}.

The dashed arcs in the figure represent the solid tori $\eta_i\times D^2$. Observe that the `outer' boundary component of $E$ is $M_L$.
Note that $E$ deformation retracts to $\overline{E}= M_L\cup (\coprod_i(\eta_i\times D^2))$, where each solid torus is attached to $M_L$ along its boundary. Hence $\overline{E}$ is obtained from $M_L$ by adding $m$ $2$-cells along the loops $\mu(\eta_i)=l_i$, and $m$ $3$-cells. Thus, by our assumption, $\phi$ extends uniquely to $\overline{\phi}:\pi_1(\overline{E})\to \G$ and hence $\ov\phi:\pi_1(E)\to\G$. Clearly the restrictions of $\overline{\phi}$ to $\pi_1(M_i)$ and $\pi_1(M_R\times \{0\})$ agree with $\phi_i$ and $\phi_R$ respectively. It follows from the third paragraph of this section that
$$
\rho(M_L,\phi) - \rho(M_R,\phi_R) = \sum^m_{i=1}\rho(M_i,\phi_i) + \sigma^{(2)}(E,\overline{\phi})- \sigma(E).
$$

Therefore most of Lemma~\ref{lem:additivity} follows from:

\begin{lem}[{\cite[Lemma 2.4]{CHL3}}]\label{lem:mickeysig} With respect to \emph{any} coefficient system, $\phi:\pi_1(E)\to \Gamma$,  the signature of the equivariant intersection form on $H_2(E;\mathbb{Z}\G)$ is zero.
\end{lem}

We want to collect, in the form of a lemma, the elementary homological properties of the cobordism $E$ that will be used in later sections.

\begin{lem}[{\cite[Lemma 2.5]{CHL3}}]\label{lem:mickeyfacts} With regard to $E$ as above, the inclusion maps induce
\begin{itemize}
\item [(1)] an epimorphism $\pi_1(M_L)\to \pi_1(E)$ whose kernel is the normal closure of the longitudes of the infecting knots $K_i$ viewed as curves $\ell_i\subset S^3-K_i\subset M_L$;
\item [(2)] isomorphisms $H_1(M_L)\to H_1(E)$ and $H_1(M_R)\to H_1(E)$;
\item [(3)] and isomorphisms $H_2(E)\cong H_2(M_L)\oplus_i H_2(M_{K_i})\cong H_2(M_R)\oplus_i H_2(M_{K_i})$.
\item [(4)] The longitudinal push-off of $\eta_i$, $\ell_{\eta_i}\subset M_L$ is isotopic in $E$ to $\eta_i\subset M_R$ and to the meridian of $K_i$, $\mu_i\subset M_{K_i}$.
\item [(5)] The longitude of $K_i$, $\ell_i\subset M_{K_i}$ is isotopic in $E$ to the reverse of the meridian of $\eta_i$, $\eta_i^{-1}\subset M_L$ and to the longitude of $K_i$ in $S^3-K_i\subset M_L$ and to the reverse of the meridian of $\eta_i$, $(\mu_{\eta_i})^{-1}\subset M_R$ (the latter bounds a disk in $M_R$).
\end{itemize}
\end{lem}

\section{Higher-Order Blanchfield forms for knots and links}\label{Blanchfieldforms}

We have seen in Lemma~\ref{lem:additivity} that an infection will have an effect on a $\rho$-invariant only if the infection circle $\eta$ survives under the map defining the coefficient system. Therefore it is important to prove \emph{injectivity} theorems concerning $\pi_1(S^3-R)\to\pi_1(B^4-\Delta)$, that is, loosely speaking, to prove that $\eta$ survives under the map
$$
j_*:\pi_1(S^3-R)^{(n)}/\pi_1(S^3-R)^{(n+1)}\to \pi_1(B^4-\Delta)^{(n)}/\pi_1(B^4-\Delta)^{(n+1)}.
$$
Higher-order Alexander modules are relevant to this task since the latter quotient can be interpreted as $H_1(W_n)$ where $W_n$ is the (solvable) covering space of $B^4-\Delta$ corresponding to the subgroup $\pi_1(B^4-\Delta)^{(n)}$. Such modules were named \emph{higher-order Alexander modules} in ~\cite{C}~\cite{COT}~\cite{Ha1}. We will employ higher-order Blanchfield linking forms on higher-order Alexander modules to find restrictions on the kernels of such maps. The logic of the technique is entirely analogous to the classical case ($n=1$): Any two curves $\eta_0, \eta_1$, say, that lie in the kernel of $j_*$ must satisfy $\mathcal{B}\ell(\eta_0,\eta_0)=\mathcal{B}\ell(\eta_0,\eta_1)=\mathcal{B}\ell(\eta_1,\eta_1)=0$ with respect to a higher order linking form $\mathcal{B}\ell$. Our major new insight is that, if the curves lie in a submanifold $S^3-K\hookrightarrow S^3-J$, a situation that arises whenever $J$ is formed from $R$ by infection using a knot $K$, then the values (above) of the higher-order Blanchfield form of $J$ can be expressed in terms of the values of the classical Blanchfield form of $K$!

Higher-order Alexander modules and higher-order linking forms for classical knot exteriors and for closed $3$-manifolds with $\beta_1(M)=1$ were introduced in ~\cite[Theorem 2.13]{COT} and further developed in ~\cite{C} and ~\cite{Lei1}. These were defined on the so called higher-order Alexander modules. Higher-order Alexander modules for \emph{links and $3$-manifolds} in general were defined and investigated in ~\cite{Ha1}. Blanchfield forms for $3$-manifolds with $\beta_1(M)>1$ were only recently defined by Leidy ~\cite{Lei3}. It is crucial to our techniques that we work with such Blanchfield forms without localizing the coefficient systems, as was investigated in ~\cite{Lei3}~\cite{Lei1}. It is in this aspect that our work deviates from that of ~\cite{CK}~\cite{COT}~\cite{COT2}. A non-localized Blanchfield form for knots also played a crucial role in ~\cite{FrT}.

First we recall that \textbf{higher-order Blanchfield linking forms} have been defined under fairly general circumstances.
\begin{thm}[{\cite[Theorem 2.3]{Lei3}}]\label{thm:blanchfieldexist} Suppose $M$ is a closed, connected, oriented $3$-manifold and $\phi:\pi_1(M)\to \Lambda$ is a PTFA coefficient system. Suppose $\mathcal{R}$ is a classical Ore localization of the Ore domain $\mathbb{Z}\Lambda$ (so $\mathbb{Z}\Lambda\subset\mathcal{R}\subset \mathcal{K}\Lambda$). Then there is a linking form:
$$
\mathcal{B}l^M_{\mathcal{R}}: TH_1(M;\mathcal{R})\to (TH_1(M;\mathcal{R}))^{\#}\equiv \overline{Hom_{\mathcal{R}}(TH_1(M;\mathcal{R}), \mathcal{K}\Lambda/\mathcal{R})}.
$$
\end{thm}

An \emph{Ore localization} of $\mathbb{Z}\Lambda$ is $\mathcal{R}=\mathbb{Z}\Lambda[S^{-1}]$ for some right-Ore set $S$ ~\cite{Ste}. When we speak of the \emph{unlocalized} Blanchfield form we mean that $\mathcal{R}=\mathbb{Z}\Lambda$ or $\mathcal{R}=\mathbb{Q}\Lambda$. $TH_1(M;\mathcal{R})$ denotes the $\mathcal{R}$-torsion submodule. In general $TH_1(M;\mathcal{R})$ need not have homological dimension one nor even be finitely-generated, and these linking forms are \emph{singular}.

Leidy analyzed the effect of an infection on the unlocalized Blanchfield forms in \cite{Lei3}~\cite{Lei1}. This generalizes the result on the classical Blanchfield form for satellite knots ~\cite{LiM}. If $L$ is obtained by infection on a link $R$ along a circle $\alpha$ using the knot $K$ and $\phi:\pi_1(M_L)\to \Lambda$ is a PTFA coefficient system, and $\mathbb{Z}\Lambda\subset\mathcal{R}\subset \mathcal{K}\Lambda$ then $\mathcal{B}l^{L}_{\mathcal{R}}$ is defined. On the other hand, by definition, exterior of the knot $K$ is a submanifold of $M_L$ and there is an induced coefficient system, that we also call $\phi$, with respect to which there is a Blanchfield linking form (first defined in ~\cite[Theorem 2.13]{COT})
$$
\mathcal{B}l_{\mathcal{R}}^K: TH_1(S^3-K;\mathcal{R})\to (TH_1(S^3-K;\mathcal{R}))^{\#}.
$$
(We note that if $\phi$ is nontrivial when restricted to $\pi_1(S^3-K)$ then $TH_1(S^3-K;\mathcal{R})=H_1(S^3-K;\mathcal{R})$. Otherwise $TH_1(S^3-K;\mathcal{R})=0$ ~\cite[Proposition 2.11]{COT}). Then it is an easy exercise for the reader using the geometric definition of these Blanchfield forms (or see ~\cite[Theorem 4.6, proof of property 1]{Lei1}), that these forms are compatible:
\begin{prop}[{\cite[Theorem 3.7]{Lei3}}]\label{prop:submanifoldcompatible} In the situation above the following diagram commutes
\begin{equation}\label{diag:compatible}
\begin{diagram}\dgARROWLENGTH=1em
\node{TH_1(S^3-K;\mathcal{R})} \arrow{e,t}{i_*}
\arrow{s,r}{\mathcal{B}l^K_\mathcal{R}}\node{TH_1(M_L;\mathcal{R})} \arrow{s,r}{\mathcal{B}l^{M_L}_\mathcal{R}}\\
\node{TH_1(S^3-K;\mathcal{R})^{\#}} \node{TH_1(M_L;\mathcal{R})^{\#}}\arrow[1]{w,t}{i^{\#}}
\end{diagram}
\end{equation}
that is, for all $x,y\in H_1(S^3-K;\mathcal{R})$
$$
\mathcal{B}l^{M_L}_{\mathcal{R}}(i_*(x),i_*(y))= \mathcal{B}l^K_\mathcal{R}(x,y).
$$
\end{prop}

Moreover, in some important situations, the induced coefficient system $\phi:\pi_1(S^3-K)\to \Lambda$ factors through, $\mathbb{Z}$, the abelianization of the knot exterior. In particular if $L$ is obtained by infection on a link $R$ along a circle $\alpha\in \pi_1(M_R)^{(k-1)}$ where $\Lambda^{(k)}=1$, then this is the case. Furthermore the higher-order Blanchfield form $\mathcal{B}l^K_\Lambda$ is merely the classical Blanchfield form on the classical Alexander module, ``tensored up.'' What is meant by this is the following. Supposing that $\phi$ is both nontrivial and factors through the abelianization, the induced map $\text{image}(\phi)\equiv\mathbb{Z}\hookrightarrow \Lambda$ is an embedding so it induces embeddings
$$
\phi:\mathbb{Q}[t,t^{-1}]\hookrightarrow \mathbb{Q}\Lambda,~~~ \phi:\mathbb{Q}(t)\hookrightarrow \mathcal{K}\Lambda,
$$
and hence an embedding (~\cite[Lemma 6.5]{CHL3})
$$
\ov\phi:\mathbb{Q}(t)/\mathbb{Q}[t,t^{-1}]\hookrightarrow \mathcal{K}\Lambda/\mathbb{Q}\Lambda.
$$

Then there is an isomorphism
$$
H_1(S^3\backslash K;\mathbb{Q}\Lambda)\cong H_1(S^3\backslash K;\mathbb{Q}[t,t^{-1}])\otimes_{\mathbb{Q}[t,t^{-1}]}\mathbb{Q}\Lambda \cong \mathcal{A}_0(K)\otimes_{\mathbb{Q}[t,t^{-1}]}\mathbb{Q}\Lambda,
$$
where $\mathcal{A}_0(K)$ is the classical (rational) Alexander module of $K$ and where $\mathbb{Q}\Lambda$ is a $\mathbb{Q}[t,t^{-1}]$-module via the map $t\to \phi(\alpha)$ ~\cite[Theorem 8.2]{C}. Moreover
$$
\mathcal{B}l_{\Lambda}^K(x\otimes 1,y\otimes 1)=\ov\phi(\mathcal{B}l^K_0(x,y))
$$
for any $x,y\in \mathcal{A}_0(K)$, where $\mathcal{B}l_0^K$ is the classical Blanchfield form on the rational Alexander module of $K$ ~\cite[Proposition 3.6]{Lei3}~\cite[Theorem 4.7]{Lei1} (see also ~\cite[Section 5.2.2]{Cha2}).

Then, finally, Leidy shows that the Blanchfield form on $M_L$ is the sum of that on $H_1(M_R)$ and that on the infecting knot $K$ (generalizing the classical result for satellites ~\cite{LiM}). We state this below although, in this paper, we shall not need this nontrivial fact that the module $H_1(M_L;\mathbb{Q}\Lambda)$ decomposes, nor even that $\mathcal{A}_0(K)\otimes_{\mathbb{Q}[t,t^{-1}]}\mathbb{Q}\Lambda$ is a submodule of it. We will only need the almost obvious fact that the inclusion of the $3$-manifolds $S^3-K_i\hookrightarrow M_L$ induces a (natural) map on the Blanchfield forms and that the induced Blanchfield form on $S^3-K$ is the classical form ``tensored up.''

\begin{thm}[{\cite[Theorem 3.7, Proposition 3.4]{Lei3}}]\label{thm:decomposition} Suppose $L=R(\alpha_i,K_i)$ is obtained by infection as above with $\alpha_i\in \pi_1(M_R)^{(k-1)}$ for all $i$. Let the zero surgeries on $R$, $L$, and $K_i$ be denoted $M_R$ $M_L$, $M_i$ respectively. Suppose $\Lambda$ is a PTFA group such that $\Lambda^{(k)}=1$. Suppose $\phi:\pi_1(M_L)\to \Lambda$ is a coefficient system. Then the inclusions induce an isomorphism
$$
H_1(M_R;S^{-1}\mathbb{Z}\Lambda)\oplus_{i\in A} H_1(S^3\backslash K_i;S^{-1}\mathbb{Z}\Lambda)\overset{i_*}{\to} H_1(M_L;S^{-1}\mathbb{Z}\Lambda).
$$
where $A=\{i~|~\phi((\alpha_i)^+\neq 1\}$. Moreover there is an isomorphism
$$
H_1(S^3\backslash K_i;\mathbb{Q}[t,t^{-1}])\otimes_{\mathbb{Q}[t,t^{-1}]}S^{-1}\mathbb{Z}\Lambda \cong H_1(S^3\backslash K_i;S^{-1}\mathbb{Z}\Lambda).
$$
Restricting to $S^{-1}\mathbb{Z}\Lambda=\mathbb{Q}\Lambda$ for simplicity, for any $x,y\in H_1(S^3\backslash K_i;\mathbb{Q}[t,t^{-1}])$,
$$
\mathcal{B}l_{\mathbb{Q}\Lambda}^{M_L}(i_*(x\otimes 1),i_*(y\otimes 1))=\ov\phi_i(\mathcal{B}l_0^i(x,y))
$$
where $\mathcal{B}l_{\Lambda}^{M_L}$ is the Blanchfield form on $M_L$ induced by $\phi$, $\mathcal{B}l_0^i$ is the classical Blanchfield form on the classical rational Alexander module of $K_i$, and
$$
\ov\phi_i: \mathbb{Q}(t)/\mathbb{Q}[t,t^{-1}]\to \mathcal{K}\Lambda/\mathbb{Q}\Lambda
$$
is the monomorphism induced by $\phi:\mathbb{Z}\to \Lambda$ sending $1$ to $\phi(\alpha_i)$.
\end{thm}

 \begin{remark}Under our hypotheses the coefficient system $\phi$ extends over the cobordism $E$, as in the discussion preceding Lemma~\ref{lem:additivity}, and there is a unique induced coefficient system $\phi_R$ on $M_R$. By Property $(4)$ of Lemma~\ref{lem:mickeyfacts}, $\alpha_i$ and its longitudinal push-off $\alpha_i^+$ are isotopic in $E$ so $\phi((\alpha_i)^+)=\phi_R(\alpha_i)$. Thus $\phi((\alpha_i)^+)\neq 1$ if and only if $\phi_R(\alpha_i\neq 1)$. Moreover, since the meridian of $K_i$ is equated to $(\alpha_i)^+$, $\phi_i(\mu_i)=\phi((\alpha_i)^+)=\phi_R(\alpha_i)$.
 \end{remark}

 The following is perhaps the key result of the paper, that we use to establish a certain  injectivity as discussed in the first paragraph of this section. Recall that the notions of $(n)$-solvable and rationally $(n)$-solvable are defined in Section~\ref{sec:appendix}. For the reader who is just concerned with proving that knots and links are not slice, replace the hypothesis below that ``$W$ is a rational $(k)$-solution for $M_L$'' with the hypothesis that ``$L$ is a slice link and $W$ is the exterior in $B^4$ of a set of slice disks for $L$.'' Such an exterior is a rational $(k)$-solution for any $k$.

\begin{thm}\label{thm:nontriviality} Suppose $L=R(\alpha_i,K_i)$ is obtained by infection. Let the zero surgeries on $R$, $L$, and $K_i$ be denoted $M_R$ $M_L$, $M_i$ respectively. Suppose $\alpha_i\in \pi_1(M_R)^{(k-1)}$ for all $i$. Suppose $W$ is a rational $(k)$-solution for $M_L$, $\Lambda$ is a PTFA group such that $\Lambda^{(k)}=1$, and $\psi:\pi_1(W)\to \Lambda$ is a nontrivial coefficient system whose restriction to $\pi_1(M_L)$ is denoted $\phi$. Let $A=\{i~|~\phi((\alpha_i)^+)\neq 1\}$. For each $i\in A$, let $P_i$ be the kernel of the composition
$$
\mathcal{A}_0(K_i)\overset{id\otimes 1}\lra  ~(\mathcal{A}_0(K_i) \otimes_{\mathbb{Q}[t,t^{-1}]}\mathbb{Q}\Lambda)\overset{i_*}{\to} H_1(M_L;\mathbb{Q}\Lambda)\overset{j_*}\to H_1(W;\mathbb{Q}\Lambda).
$$
Then $P_i\subset P_i^\perp$ with respect to $\mathcal{B}l_0^i$, the classical Blanchfield linking form on the rational Alexander module, $\mathcal{A}_0(K_i)$, of $K_i$.
\end{thm}

\begin{remark} Under the hypotheses of Theorem~\ref{thm:nontriviality}, the coefficient system extends over the cobordism $E$ of Figure~\ref{fig:mickey} and hence extends to $\pi_1(M_R)$. If this extension is (sloppily) also called $\phi$ then $\phi(\alpha_i)=\phi((\alpha_i)^+)$ since $\alpha_i$ and its longitude $(\alpha_i)^+$ are isotopic in $M_R$ and hence freely homotopic in $E$.
\end{remark}

\begin{proof}[Proof of Theorem~\ref{thm:nontriviality}] We need:

\begin{lem}\label{lem:fourmanBlanch} There is a Blanchfield form, $\mathcal{B}l^{rel}$,
$$
\mathcal{B}l^{rel}:~TH_2(W,\partial W;\mathcal{R})\to TH_1(W)^{\#}
$$
such that the following diagram, with coefficients in $\mathcal{R}$ unless specified otherwise, is commutative up to sign:
%\begin{equation}\label{diag:compatible2}
$$
\begin{diagram}\label{diag:compatible2}\dgARROWLENGTH=1em
\node{TH_2(W,\partial W;\mathcal{R})} \arrow{e,t}{\partial_*}
\arrow{s,r}{\mathcal{B}l^{rel}_\mathcal{R}}\node{TH_1(M;\mathcal{R})} \arrow{s,r}{\mathcal{B}l^{M}_\mathcal{R}}\\
\node{TH_1(W;\mathcal{R})^{\#}} \arrow{e,t}{j^{\#}}\node{TH_1(M;\mathcal{R})^{\#}}
\end{diagram}
$$
%\end{equation}
\end{lem}

\begin{proof}[Proof of Lemma~\ref{lem:fourmanBlanch}] (See also ~\cite[Lemmas 3.2, 3.3]{Cha3})  Consider the following commutative diagram where homology and cohomology is with $\mathcal{R}$ coefficients unless specified and $\mathcal{K}$ denotes the quotient field of $\mathcal{R}$:
$$
\begin{diagram}[small]
\node{H_3(W,M;\mathcal{K})} \arrow[2]{e,t}{\partial_*} \arrow[2]{s,l}{P.D.} \arrow{se} \node[2]{H_2(M;\mathcal{K})} \arrow{s,-} \arrow{se}\\
\node[2]{H_3(W,M;\mathcal{K}/\mathcal{R})} \arrow[2]{e} \arrow[2]{s} \node{} \arrow{s} \node{H_2(M;\mathcal{K}/\mathcal{R})} \arrow[2]{s} \\
\node{\overline{H^1(W;\mathcal{K})}} \arrow[2]{s,l}{\kappa} \arrow{se} \arrow{e,-} \node{} \arrow{e} \node{\overline{H^1(M;\mathcal{K})}} \arrow{s,-} \arrow{se} \\
\node[2]{\overline{H^1(W;\mathcal{K}/\mathcal{R})}} \arrow[2]{e} \arrow[2]{s} \node{} \arrow{s} \node{\overline{H^1(M;\mathcal{K}/\mathcal{R})}} \arrow[2]{s} \\
\node{\overline{\Hom_{\mathcal{R}}(H_1(W),\mathcal{K})}} \arrow[2]{s,l}{\iota} \arrow{se} \arrow{se} \arrow{e,-} \node{} \arrow{e} \node{\overline{\Hom_{\mathcal{R}}(H_1(M),\mathcal{K})}} \arrow{s,-} \arrow{se} \\
\node[2]{\overline{\Hom_{\mathcal{R}}(H_1(W),\mathcal{K}/\mathcal{R})}} \arrow[2]{e} \arrow[2]{s} \node{} \arrow{s} \node{\overline{\Hom_{\mathcal{R}}(H_1(M),\mathcal{K}/\mathcal{R})}} \arrow[2]{s} \\
\node{\overline{\Hom_{\mathcal{R}}(TH_1(W),\mathcal{K})}} \arrow{se} \arrow{se} \arrow{e,-} \node{} \arrow{e} \node{\overline{\Hom_{\mathcal{R}}(TH_1(M),\mathcal{K})}} \arrow{se} \\
\node[2]{\overline{\Hom_{\mathcal{R}}(TH_1(W),\mathcal{K}/\mathcal{R})}} \arrow[2]{e,t}{j^{\#}} \node[2]{\overline{\Hom_{\mathcal{R}}(TH_1(M),\mathcal{K}/\mathcal{R})}} \\
\end{diagram}
$$
%$$
%\begin{diagram}
%\node{H_3(W,M;\mathcal{K})} \arrow[2]{e,t}{\partial_*} \arrow[2]{s,l}{\iota \circ \kappa \circ P.D.} \arrow{se} %\node[2]{H_2(M;\mathcal{K})} \arrow{s,-} \arrow{se} \\
%\node[2]{H_3(W,M;\mathcal{K}/\mathcal{R})} \arrow[2]{e} \arrow[2]{s} \node{} \arrow{s} %\node{H_2(M;\mathcal{K}/\mathcal{R})} \arrow[2]{s} \\
%\node{\overline{\Hom_{\mathcal{R}}(TH_1(W),\mathcal{K})}} \arrow{se} \arrow{e,-} \node{} \arrow{e} %\node{\overline{\Hom_{\mathcal{R}}(TH_1(M),\mathcal{K})}} \arrow{se} \\
%\node[2]{\overline{\Hom_{\mathcal{R}}(TH_1(W),\mathcal{K}/\mathcal{R})}} \arrow[2]{e,t}{j^{\#}} %\node[2]{\overline{\Hom_{\mathcal{R}}(TH_1(M),\mathcal{K}/\mathcal{R})}} \\
%\end{diagram}
%$$
where $\iota$ is the map induced from the inclusion map of the torsion submodule. Since
$$
\overline{\Hom_{\mathcal{R}}(TH_1(W;\mathcal{R}),\mathcal{K})}=0,
$$
it follows that the image of $H_3(W,M;\mathcal{K}) \to H_3(W,M;\mathcal{K}/\mathcal{R})$ is contained in the kernel of the composition $\iota \circ \kappa\circ P.D.$. Furthermore, from the exact sequence,
$$
H_3(W,M;\mathcal{K})\overset{\pi}{\to} H_3(W,M;\mathcal{K}/\mathcal{R})\to H_2(W,M;\mathcal{R})\to H_2(W,M;\mathcal{K})
$$
since $H_2(W,M;\mathcal{K})$ is $\mathcal{R}$-torsion-free, $TH_2(W,M;\mathcal{R})$ is isomorphic to the cokernel of $\pi$. It follows that there is a well-defined map $\mathcal{B}l^{rel}_\mathcal{R}:
TH_2(W,M;\mathcal{R}) \to TH_1(W;\mathcal{R})^{\#}$. Similarly, since
$$
\overline{\Hom_{\mathcal{R}}(TH_1(M;\mathcal{R}),\mathcal{K})}=0,
$$
there is a well-defined map
$\mathcal{B}l^{M}_\mathcal{R}: TH_1(M,\mathcal{R}) \to TH_1(M;\mathcal{R})^{\#}$ such that the following diagram commutes.
%$$
%\begin{diagram}[small]
%\node{H_3(W,M;\mathcal{K}/\mathcal{R})} \arrow{s,l}{P.D.} \arrow[2]{e} \arrow{sse} %\node[2]{H_2(M;\mathcal{K}/\mathcal{R})} \arrow{s} \arrow{sse} \\
%\node{\overline{H^1(W;\mathcal{K}/\mathcal{R})}} \arrow[2]{s,l}{\kappa} \arrow[2]{e} %\node[2]{\overline{H^1(M;\mathcal{K}/\mathcal{R})}} \arrow{s,-} \\
%\node[2]{TH_2(W,M)} \arrow{ssw,r}{\mathcal{B}l^{rel}_\mathcal{R}} \arrow[2]{e,t}{\partial_*\hspace{.25in}} \node{} %\arrow{s} \node{TH_1(M)} \arrow{ssw,r}{\mathcal{B}l^{M}_\mathcal{R}}\\
%\node{H_1(W)^{\#}} \arrow{s,l}{\iota} \arrow[2]{e} \node[2]{H_1(M)^{\#}} \arrow{s} \\
%\node{TH_1(W)^{\#}} \arrow[2]{e,t}{j^{\#}} \node[2]{TH_1(M)^{\#}} \\
%\end{diagram}
%$$
$$
\begin{diagram}[small]
\node{H_3(W,M;\mathcal{K}/\mathcal{R})} \arrow[2]{s,l}{\iota \circ \kappa
\circ P.D.} \arrow[2]{e} \arrow{se} \node[2]{H_2(M;\mathcal{K}/\mathcal{R})} \arrow{s,-} \arrow{se} \\
\node[2]{TH_2(W,M)} \arrow{sw,r}{\mathcal{B}l^{rel}_\mathcal{R}} \arrow[2]{e,t}{\partial_*\hspace{.25in}} \node{} \arrow{s} \node{TH_1(M)} \arrow{sw,r}{\mathcal{B}l^{M}_\mathcal{R}}\\
\node{TH_1(W)^{\#}} \arrow[2]{e,t}{j^{\#}} \node[2]{TH_1(M)^{\#}} 
\end{diagram}
$$
\end{proof}

The following result was proved in ~\cite[Lemma 4.5, Theorem 4.4]{COT} in the special case that $\beta_1(M)=1$. It was proved in more generality in ~\cite[Theorem 6.3]{CHL3} except that there the proof of Lemma~\ref{lem:fourmanBlanch} was omitted.

\begin{lem}\label{lem:selfannihil} Suppose $M$ is connected and is rationally
$(k)$-solvable via $W$ and
$\phi:\pi_1(W)\ra\Lambda$ is a non-trivial coefficient system where
$\Lambda$ is a PTFA group with $\Lambda^{(k)}=1$. Let $\mathcal{R}$ be an Ore localization of $\mathbb{Z}\Lambda$ so $\mathbb{Z}\Lambda\subset\mathcal{R}\subset \mathcal{K}\Lambda$. Then
$$
TH_2(W,M;\mathcal{R})\xrightarrow{\partial}TH_1(M;\mathcal{R})\xrightarrow{j_{\ast}} TH_1(W;\mathcal{R})
$$
is exact. Moreover, any submodule $P\subset \text{kernel}~j_*$ satisfies
$P\subset ($ker~$j_*)^\perp\subset P^\perp$ with respect to the Blanchfield form on $TH_1(M;\mathcal{R})$.
\end{lem}

We can now finish the proof of Theorem~\ref{thm:nontriviality}. Suppose $x,y\in P_i$ as in the statement. Let $\mathcal{R}=\mathbb{Q}\Lambda$, $M=M_L$ and let $P$ be the submodule of $H_1(M_L;\mathbb{Q}\Lambda)$  generated by $\{i_*(x\otimes 1),i_*(y\otimes 1)\}$. Then $P\subset \text{kernel}~j_*$. Apply Lemma~\ref{lem:selfannihil} to conclude that
$$
\mathcal{B}l^{M_L}_{\mathbb{Q}\Lambda}(i_*(x\otimes 1)),(i_*(y\otimes 1))=0.
$$
By Theorem~\ref{thm:decomposition},
$$
\ov\phi_i(\mathcal{B}l_o^i(x,y))=0.
$$
Since $\ov\phi$ is a monomorphism, it follows that $\mathcal{B}l_o^i(x,y)=0$. Thus $P_i\subset P_i^\perp$ with respect to the classical Blanchfield form on $K_i$. This concludes the proof of Theorem~\ref{thm:nontriviality}.
\end{proof}

\section{Iterated Bing doubles and first-order $L^{(2)}$-signatures}\label{sec:Bingdoubles}

In this section we investigate higher-order signature invariants that obstruct any iterated Bing double of $K$ from being a topologically slice link. We first state and prove the simplest results and later generalize.

Harvey and Cha-Livingston-Ruberman showed that classical signatures of $K$, which we call $0^{th}$-order signatures, obstruct $BD(K)$ from being a slice link. These signatures vanish if $K$ is an algebraically slice knot. Here we show that certain higher-order signatures of $K$, similar to Casson-Gordon invariants, that we call \emph{first-order signatures} of $K$,  obstruct $BD(K)$ from being a slice link. These were first defined in ~\cite{CHL1A} (see also ~\cite{CHL3}). To define these, suppose $K$ is an oriented knot, let $G=\pi_1(M_K)$ and let $\mathcal{A}_0=\mathcal{A}_0(K)$ be its classical rational Alexander module. Note that since the longitudes of $K$ lie in $\pi_1(S^3-K)^{(2)}$,
$$
\mathcal{A}_0\equiv G^{(1)}/G^{(2)}\otimes_{\mathbb{Z}[t,t^{-1}]}\mathbb{Q}[t,t^{-1}].
$$
Each submodule $P\subset \mathcal{A}_0$ corresponds to a unique metabelian quotient of $G$,
$$
\phi_P:G\to G/\tilde{P},
$$
by setting
$$
\tilde{P}\equiv \{x~| x\in \text{kernel}(G^{(1)}\to G^{(1)}/G^{(2)}\to \mathcal{A}_0/P)\}.
$$
\noindent Note that $G^{(2)}\subset \tilde{P}$ so $G/\tilde{P}$ is metabelian. Therefore to any such submodule $P$ there corresponds a real number, the Cheeger-Gromov invariant, $\rho(M_K, \phi_P:G\to G/\tilde{P})$.

\begin{defn}\label{defn:highordersignatures} The \textbf{first-order $\boldsymbol{L^{(2)}}$-signatures} of a knot $K$ are the real numbers
$\rho(M_K, \phi_P)$ where $P\subset \mathcal{A}_0(K)$ satisfies $P\subset P^\perp$.
\end{defn}

Since $P=0$ (the case $\tilde{P}=G^{(2)}$) always satisfies $P\subset P^\perp$, we give a special name to the first-order signature corresponding to this case.

\begin{defn}\label{defn:rho1} $\boldsymbol{\rho^1(K)}$ of a knot $K$ is the first-order $L^{(2)}$-signature given by the Cheeger-Gromov invariant $\rho(M_K, \phi:G\to G/G^{(2)})$.
\end{defn}

\begin{figure}[htbp]
\setlength{\unitlength}{1pt}
\begin{center}
\begin{picture}(143,151)
\put(0,0){\includegraphics{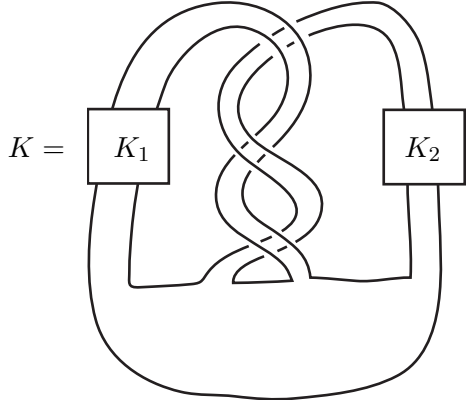}}
\put(-30,92){$K=$}
\put(10,92){$K_1$}
\put(120,92){$K_2$}
\end{picture}
\end{center}
\caption{A genus $1$ algebraically slice knot $K$}\label{fig:examplehighersigs}
\end{figure}

\begin{ex}\label{ex:first-ordersigs} Consider the knot $K$ in Figure~\ref{fig:examplehighersigs}. This knot is obtained from the ribbon knot $R=9_{46}$ by two infections on the band meridians $\alpha, \beta$ (as in the left-hand side of Figure~\ref{fig:Rdoubling}). Thus $\{\alpha, \beta \}$ is a basis of $\mathcal{A}_0(K)= \mathcal{A}_0(9_{46})$. There are $3$ submodules $P$ for which $P\subset P^\perp$, namely $P_0=0$, $P_1=\left<\alpha\right>$ and $P_2=\left<\beta\right>$. We may apply Lemma~\ref{lem:additivity} to show
$$
\rho(M_K, \phi_P)=\rho(M_{R}, \phi_P)+\epsilon^1 _P\rho_0(K_1)+\epsilon^2 _P\rho_0(K_2)
$$
where $\epsilon^1_P$ is $0$ or $1$ according as $\phi_P(\alpha)=1$ or not (similarly for $\epsilon^2_P$). For our example $\phi_{P_1}(\alpha)=1$ and $\phi_{P_1}(\beta)\neq 1$. Similarly $\phi_{P_2}(\beta)=1$ and $\phi_{P_2}(\alpha)\neq 1$. By contrast $\phi_{P_0}(\alpha)\neq 1$ and $\phi_{P_0}(\beta)\neq 1$. Moreover $P_1$ corresponds to the kernel $\tilde{P}_1$, of $\pi_1(S^3-R)\to \pi_1(B^4-\Delta_1)/\pi_1(B^4-\Delta_1)^{(2)}$ for the ribbon disk $\Delta_1$ for $R$ obtained by ``cutting the $\alpha$-band.'' (Similarly for $P_2$.) Thus in both cases the maps $\phi_P$ on $M_{R}$ extend over ribbon disk exteriors. Consequently $\rho(M_{R}, \phi_P)=0$ for $P=P_1$ and $P=P_\beta$, by Theorem \ref{thm:oldsliceobstr}. Of course $\rho(M_R, \phi_{P_0})=\rho^1(9_{46})$ by definition. Putting this all together we see that the first-order signatures of the knot $K$ are $\{\rho^1(9_{46})+\rho_0(K_1)+\rho_0(K_2),\rho_0(K_2),\rho_0(K_1)\}$.
\end{ex}

\begin{ex}\label{ex:first-ordersigs89}  Consider the ribbon knot $8_9$, pictured on the left-hand side of Figure $4.2$ ~\cite{Lam}. A ribbon move is indicated by the small dotted arc. We will show that all of its first-order signatures are zero.
The $8_9$ knot is a ribbon knot and is fully ($\pm$) amphicheiral with Alexander polynomial $p(t)q(t)$, where $p(t)=t^3-2t^2+t-1$ and $q(t)=t^3-t^2+2t-1$ ~\cite[p.263, 270, 279]{Ka3}. Since $p$ and $q$ are irreducible and distinct (up to units) in the PID $\mathbb{Q}[t,t^{-1}]$, the Alexander module of $8_9$ is cyclic of order $pq$. It follows that $\mathcal{A}_0(8_9)$ has precisely $3$ proper submodules: $P_0=0, P_1=<p>, $ and $P_2=<q>$ and hence $3$ first-order signatures. The first-order signature corresponding to $P_0$ (what we  call $\rho^1(8_9)$) is zero by the following result of the authors.

\begin{figure}[htbp]
\setlength{\unitlength}{1pt}
\begin{center}
\begin{picture}(316,135)
\put(0,0){\includegraphics{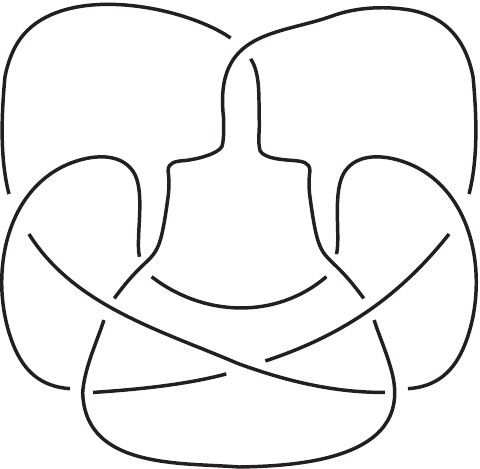}}
\put(188,0){\includegraphics{8_9prettierboxes.pdf}}
\put(-30,80){$8_9~=$}
\put(57,130){............}
\put(251,102){$K_1$}
\put(225,72){$K_1$}
\put(160,80){$K=$}
\end{picture}
\end{center}
\caption{}\label{fig:8_9}
\end{figure}

\begin{prop}[{\cite[Proposition 5.3]{CHL1A}\cite[Proposition 3.4]{CHL3}}]\label{prop:amphi} If a $3$-manifold $M$ admits an orientation-reversing homeomorphism $h$, then $\rho(M,\phi)=0$ for any $\phi$ whose kernel is invariant under $h_*$, in particular if the kernel is a characteristic subgroup.
\end{prop}

\noindent Moreover, since $8_9$ is a ribbon knot, it admits a slice disk, $\Delta_1$. It is a classical result that the kernel of
$$
(i_1)_*:\mathcal{A}_0(8_9)\to H_1(B^4-\Delta_1;\mathbb{Q}[t,t^{-1}])
$$
is self-annihilating with respect to the Blanchfield form, so this kernel is either $P_1$ or $P_2$, say $P_1$ for specificity. It follows that the kernel of the inclusion-induced map
$$
(i_1)_*:G\to G/G^{(2)}\to\pi_1(B^4-\Delta_1)/\pi_1(B^4-\Delta_1)^{(2)}_r
$$
is $\tilde{P}_1$. Since $\phi_{P_1}$ extends to the exterior of this slice disk the first-order signature corresponding to $P_1$ is zero. This leaves only $P_2$. Consider a homeomorphism $\tilde{h}$ of $B^4$ that restricts to a reflection $h$ on $S^3$. The image $h(8_9)$ is the mirror image, $\overline{8}_9$, of $8_9$ and the image of $\Delta_1$ is a ribbon disk, $\Delta_2$ for $\overline{8}_9$. Since $8_9$ is isotopic to its mirror image, this can be viewed as another ribbon disk for $8_9$. The kernel of the map $(i_2)_*$, as above, is $h_*(\tilde{P_1})$ and must be either $\tilde{P_1}$ or $\tilde{P_2}$. If it is $\tilde{P}_2$ then the first-order signature corresponding to $P_2$ vanishes, since $\phi_{P_2}$ extends to the exterior of this slice disk. If not, then $h_*(\tilde{P_1})=\tilde{P_1}$ and consequently $h_*(\tilde{P_2})=\tilde{P_2}$. Then, again since $K$ is amphichiral, and $h$ preserves the kernel of $\phi_{P_2}$, the first-order corresponding to $P_2$ vanishes by Proposition~\ref{prop:amphi}.
Thus all of the first-order signatures of $8_9$ are zero. (Note: in fact it can be shown that $h_*(\tilde{P_1})=\tilde{P_2}$.)
\end{ex}

%\begin{ex}\label{ex:0first-ordersigs89T} Consider the ribbon knot $R$ on the right-hand side of Figure?? %where $T$ is the trefoil knot or any other knot such that $\rho_0(T)\neq 0$. We claim that %$\rho^1(R)=\rho_0(T)\neq 0$. For since $R$ is obtained from $8_9$ by a single infection a ribbon-band %meridian, which represents a non-trivial element in the Alexander module of $8_9$, we may apply %Lemma~\ref{lem:additivity} to yield
%$$
%\rho^1(R)=\rho(M_R, \phi_0)=\rho^1(8_9)+\rho_0(T)=\rho_0(T).
%$$
%\end{ex}

\begin{ex}\label{ex:nonzerofirst-ordersigs} Consider the family of algebraically slice knots, $K$, shown on the right-hand side of Figure~\ref{fig:8_9}. Suppose $\rho_0(K_1)\neq 0$. Then we claim that all of the first-order signatures of $K$ are \emph{non-zero}. Since $K$ has the same Alexander module as $8_9$ it has $3$ first-order signatures. First note that $K$ is obtained from $8_9$ by two infections. The infection using the upper copy of $K_1$ is done along a curve that represents a generator of the cyclic module $\mathcal{A}_0(8_9)$. Such a curve cannot lie in \emph{any} submodule $P\subset P^\perp$. The infection using the lower copy of $K_1$ is done along a generator of $P_1$, hence does not lie in $P_2$. Since all of the first-order signatures of $8_9$ are zero, by Lemma~\ref{lem:additivity} the first-order signatures of $K$ corresponding to $\{P_0,P_1,P_2\}$ are, respectively, $\{2\rho_0(K_1),\rho_0(K_1),2\rho_0(K_1)\}$, each of which is non-zero.
\end{ex}

We will now show that the first-order signatures of an arbitrary knot $K$, like the ordinary signatures, obstruct each \emph{iterated Bing double} of $K$ from being a (topologically) slice link. This improves on Harvey's theorem which showed this same fact for the integral of the classical signatures ~\cite[Corollary 5.6]{Ha2}. There are several ways to define \textbf{iterated Bing Doubling}. In the most general way, one doubles one component at a time. However for simplicity, let us focus on the notion of Bing doubling wherein we double \emph{every} component. Then the $n$-fold iterated Bing double of $K$, $BD^{n}(K)$, is a $2^n$ component link. Note that once we show that none of these restricted Bing doubles is slice then it follows that none of the more general iterated Bing doubles is slice.

\begin{thm}\label{thm:Bingdouble} Let $K$ be an arbitrary knot. If the $n$-fold iterated Bing double of $K$ ($n\geq 1$) is topologically slice in a rational homology $4$-ball (or more generally is a rationally $(n+1.5)$-solvable link) then one of the first-order signatures of $K$ is zero.
\end{thm}

\noindent Before proving Theorem~\ref{thm:Bingdouble}, we establish two corollaries.

\begin{figure}[htbp]
\setlength{\unitlength}{1pt}
\begin{center}
\begin{picture}(159,101)
\put(0,0){\includegraphics{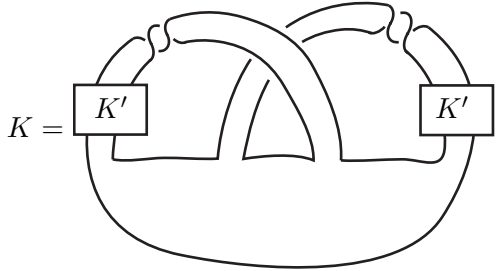}}
\put(-25,50){$K=$}
\put(8,57){$K'$}
\put(137,57){$K'$}
%\put(-40,70){$K=$}
%\put(7,95){$K'$}
%\put(145,94){$K'$}
\end{picture}
\end{center}
\caption{Order 2 in algebraic concordance group}\label{fig:examplehighersigsfigeight}
\end{figure}

\begin{cor}\label{cor:Bingdouble2} If $K$ is one of the algebraically slice knots of Example~\ref{ex:nonzerofirst-ordersigs} then the $n$-fold iterated Bing double of $K$ is $(n+1)$-solvable (requires also that Arf$(K_1)=0$) but not slice nor even rationally $(n+1.5)$-solvable. Similarly, if $K$ is the knot of Figure~\ref{fig:examplehighersigsfigeight} where $\rho_0(K')\neq 0$ then no iterated Bing double of $K$ is topologically slice (nor even rationally $(n+1.5)$-solvable).
\end{cor}

\begin{proof}
%\begin{proof} [Proof of Corollary~\ref{cor:Bingdouble2}]
Let $R$ denote the $8_9$ knot. Then our knot $K$ is obtained from the zero solvable knot $K_1$ by applying the $R$-operator along two curves representing elements of the commutator subgroup. Hence $K$ is $(1)$-solvable by Lemma~\ref{lem:nsolv}. Then $BD^n(K)$  is obtained by an infection on the trivial link (using $K$) along a curve in the $n^{th}$-derived (see the proof of Theorem~\ref{thm:Bingdouble}) and so $BD^n(K)$ is an $(n+1)$-solvable link by Lemma~\ref{lem:nsolv}. Apply Theorem~\ref{thm:Bingdouble} to conclude that, if $BD^n(K)$  were slice or even rationally $(n+1.5)$-solvable, then one of the first-order signatures of $K$ would vanish. The result now follows immediately from Example~\ref{ex:nonzerofirst-ordersigs}.

 For the knot of Figure~\ref{fig:examplehighersigsfigeight} there is only one submodule $P\subset P^{\perp}$ for the Alexander module of the figure-eight knot, namely $P=0$. Therefore there is only one first-order signature for the pictured knot $K$, namely $\rho^1(\text{figure-eight})+2\rho_0(K')$. Since the figure-eight knot is amphichiral, $\rho^1(\text{figure-eight})=0$, so this first-order signature is non-zero.
\end{proof}

\begin{cor}\label{cor:Bingdouble}
If $K$ is one of the algebraically slice knots of Figure~\ref{fig:examplehighersigs}, with $K_1=K_2$, and some iterated Bing double of $K$ is a slice link (or even $(n+1.5)$-solvable) then $\rho_o(K_1)\in \{0,(-1/2)\rho^1(9_{46})\}$. Therefore if $\rho_0(K_1)\notin \{0,(-1/2)\rho^1(9_{46})\}$ and Arf($K_1$)$=0$, then the $n$-fold iterated Bing double of $K$ is $(n+1)$-solvable but not slice nor even rationally $(n+1.5)$-solvable.
\end{cor}

\begin{proof}
%\begin{proof} [Proof of Corollary~\ref{cor:Bingdouble}]
Apply Theorem~\ref{thm:Bingdouble} to the knot of Figure~\ref{fig:examplehighersigs}, with $K_1=K_2$ to conclude that one of the first-order signatures of $K$ must vanish. By Example~\ref{ex:first-ordersigs}, the first-order signatures of $K$ are $\{\rho_0(K_1),2\rho_0(K_1)+\rho^1(9_{46})\}$. The claimed results follow.
\end{proof}

\begin{proof}[Proof of Theorem~\ref{thm:Bingdouble}] Let $L=BD^{n}(K)$ for some $n\geq 1$ and $M=M_{L}$. Suppose $M$ is rationally $(n+1.5)$-solvable via $V$. We shall show that one of the first-order signatures of $K$ is zero.

Recall that $BD(K)$ can be obtained from the trivial link of two components by infection on the circle $\alpha$ shown dashed in Figure~\ref{fig:bingeta}, using $K$ as the infecting knot. This curve $\alpha$ can be expressed as $[x,y]$ in the fundamental group of the zero surgery on the trivial link where $x$ and $y$ are the meridians. If one now doubles each component of this trivial link, then the image of the curve $\alpha$ becomes a curve that represents the double commutator $[[x,x'],[y,y']]$ for suitably chosen meridians. Continuing in this manner, one sees that the iterated Bing double $L$ can be obtained from the trivial $2^n$ component link $T$ by a single infection, using the knot $K$, along a circle $\alpha$ representing, in $\pi_1(M_T)$, an element in $F^{(n)}$ but not in $F^{(n+1)}$. At this point we note that we need not assume that we are dealing with an iterated Bing double, but rather this previous sentence is all that we need assume. Thus our proof is really going to prove:

\begin{thm}\label{thm:basiclink} Suppose $T$ is a trivial link of $m$ components, $n\geq 1$ and $\alpha$ is an unknotted circle in $S^3-T$ that represents an element in $F^{(n)}-F^{(n+1)}$ where $F=\pi_1(S^3-T)$, and $L$ denotes $T(\alpha,K)$, the result of infection of $T$ along $\alpha$ using the knot $K$. If $L$ is topologically slice in a rational homology $4$-ball (or is even a rationally $(n+1.5)$-solvable link) then one of the first order signatures of $K$ is zero.
\end{thm}

Proceeding with the proof of Theorem~\ref{thm:basiclink} and hence of Theorem~\ref{thm:Bingdouble}, since $L=T(\alpha,K)$, there exists a cobordism $E$ as in Figure~\ref{fig:mickey} whose boundary is $M_T\sqcup M_K\sqcup -M$. We form a $4$-manifold $W$ as follows. Cap off $M\subset \partial E$ using $V$. Thus $\partial W=M_K\cup M_T$. Let $\pi=\pi_1(W)$ and consider $\phi:\pi\to\pi/\pi^{(n+2)}_r$. In the case that $V$ is a slice disk exterior then we can apply Theorem~\ref{thm:linksliceobstr} to conclude that
$$
\rho(M,\phi)=0.
$$
If $V$ is merely a rational $(n+1.5)$-solution, we would like to apply Theorem~\ref{thm:rho=0} to arrive at the same conclusion. But we must first verify that $L$ satisfies the conditions of Lemma~\ref{lem:rank}. This requires only that $\phi(\ell_K)=1$. This is certainly the case since, by property $(5)$ of Lemma~\ref{lem:mickeyfacts}, $\ell_K$ is identified with the reverse of meridian of $\alpha$ which bounds a disk in $M_T$, hence is null-homotopic in $W$. Let $\ov\phi$ be restriction of $\phi$ to $\pi_1(M_K)$ and $\phi_T$ denote the restriction of $\phi$ to $\pi_1(M_T)$. Thus, by Lemma~\ref{lem:additivity}
$$
\rho(M_K,\ov\phi)+ \rho(M_T,\phi_T)=0.
$$
Since $T$ is a trivial link, $M_T=\partial Y$ where $Y$ is a boundary connected-sum of copies of $S^1\times B^3$. Since $\pi_1(\partial Y)\cong \pi_1(Y)$, $\phi_T$ extends to $Y$. Hence by Theorem~\ref{thm:linksliceobstr},
$$
\rho(M_T,\phi_T)=0.
$$
Therefore
$$
\rho(M_K,\ov\phi)=0.
$$

It remains only to identify $\rho(M_K,\ov\phi)$ as one of the first-order signatures of $K$. First note that the meridian of $K$ is isotopic in $E$ to the infection circle $\alpha$ in $M_T$. Since $\alpha\in \pi_1(S^3-T)^{(n)}$, this meridian represents an element of $\pi_1(E)^{(n)}$ and hence an element of $\pi^{(n)}$. Since $G\equiv\pi_1(M_K)$ is normally generated by this meridian,
$$
i_*(G)\subset \pi^{(n)}
$$
and so
$$
i_*(G^{(2)})\subset \pi^{(n+2)}.
$$
Consequently $\ov\phi$ factors through $G/G^{(2)}$ and the image of $\ov\phi$ is contained in $\pi^{(n)}/\pi^{(n+2)}_r$. By Property $2$ of Proposition~\ref{prop:rho invariants}, $\rho(M_K,\ov\phi)$ depends only on the image of $\ov\phi$. Thus
$$
\rho(M_K,\ov\phi)=\rho(M_K,G\to G/G^{(2)}\to G/\tilde{P})
$$
where $\tilde{P}=\text{ker}\ov\phi$. Therefore we need only characterize $\tilde{P}$. To this end, let $\tilde{\pi}=\pi_1(V)$. From property $(1)$ of Lemma~\ref{lem:mickeyfacts}
$$
\pi_1(M)\to \pi_1(E)
$$
is surjective with kernel the normal closure of the longitude $\ell_K$ of $K$ (here we are considering that $S^3-K\subset M$). Therefore the kernel of the map
$$
\tilde{\pi}\to \pi
$$
induced by the inclusion $V\hookrightarrow V\cup E$ is the normal closure of $\ell_K$. We claim that this induces an isomorphism
$$
\tilde{\pi}/\tilde{\pi}^{(n+2)}_r\cong \pi/\pi^{(n+2)}_r.
$$
This will follow if we show $\ell_K\in \tilde{\pi}^{(n+2)}$.  Recall that $\alpha\in \pi_1(S^3-T)^{(n)}$. It follows, as shown in ~\cite[Proof of Theorem 8.1]{C} that a stronger fact holds, namely that the longitudinal push-off of $\alpha$, $\ell_\alpha$, lies in $\pi_1(M)^{(n)}$. But $\ell_\alpha$ is identified to the meridian, $\mu_K$, of $S^3-K\subset M$. Since $\ell_K\in \pi_1(S^3-K)^{(2)}$ and $\pi_1(S^3-K)$ is normally generated by $\mu_K$,
$$
\ell_K\in \pi_1(M)^{(n+2)}\subset \tilde{\pi}^{(n+2)},
$$
as required.
Hence
$$
\tilde{P}=\text{ker}(G\to \tilde{\pi}/\tilde{\pi}^{(n+2)}_r).
$$
Moreover, since the copy of $S^3-K$ that is a subset of $M_K$ and the copy of $S^3-K$ that is a subset of $M$ are isotopic in $E$, we are now free to think of $G$ as $\pi_1$ of the latter copy (modulo the longitude).

Now consider $\Lambda=\tilde{\pi}/\tilde{\pi}^{(n+1)}_r\cong \pi/\pi^{(n+1)}_r$ and $\psi:\tilde{\pi}\to \Lambda$. We seek to apply Theorem~\ref{thm:nontriviality} to $L=T(\alpha,K)$, $\alpha\in \pi_1(S^3-T)^{(n)}$, $k=n+1$ and the rational $(n+1)$-solution $V$ for $M$. To apply Theorem~\ref{thm:nontriviality}, we first need to verify that $\psi(\alpha)\neq 1$.

Consider the inclusion $i:M_T\to W$. By property $(2)$ of Lemma~\ref{lem:mickeyfacts} and since $V$ is a rational $(n)$-solution, this map induces an isomorphism on $H_1(-;\mathbb{Q})$. By property $(3)$ of Lemma~\ref{lem:mickeyfacts}
$$
H_2(W;\mathbb{Q})\cong H_2(V;\mathbb{Q})\oplus i_*(H_2(M_K;\mathbb{Q}).
$$
Since $V$ is a rational $(n)$-solution, $H_2(V;\mathbb{Q})$ has a basis consisting of surfaces $\Sigma$ wherein $\pi_1(\Sigma)\subset \pi^{(n)}$. $H_2(M_K)$ is represented by a capped off Seifert surface $\overline{\Sigma}$ for $K$. Since $\pi_1(M_K)$ is normally generated by the meridian of $K$, which lies in $\pi^{(n)}$, $\pi_1(\overline{\Sigma})\subset \pi^{(n)}$. Thus, by ~\cite[Theorem 2.1]{CH2}, there is a monomorphism
$$
i_H:\pi_1(M_T)/\pi_1(M_T)^{(n+1)}_H\hookrightarrow \pi/\pi^{(n+1)}_H
$$
where the subscript $H$ denotes Harvey's torsion-free derived series ~\cite[Section 2]{Ha2}. Since the rational derived series is contained in the torsion-free derived series we have the commutative diagram
\begin{equation}\label{diag:harvey}
\begin{diagram}\dgARROWLENGTH=1em
\node{\pi_1(M_T)/\pi_1(M_T)^{(n+1)}_r} \arrow{e,t}{i_*}
\arrow{s,l}{\pi}\node{\pi/\pi^{(n+1)}_r} \arrow{s,l}{}\arrow{e,t}{\cong}\node{\Lambda}\\
\node{\pi_1(M_T)/\pi_1(M_T)^{(n+1)}_H}\arrow{e,t}{i_H} \node{\pi/\pi^{(n+1)}_H}
\end{diagram}
\end{equation}
From this diagram we see that if $\alpha\in \pi_1(M_T)$ mapped to zero in $\Lambda$ then
$\pi(\alpha)=1$ meaning that $\alpha\in \pi_1(M_T)^{(n+1)}_H$. But this contradicts our hypothesis on $\alpha$ since, for the free group $\pi_1(M_T)$, the torsion-free derived series coincides with the derived series ~\cite[Proposition 2.3]{Ha2}. Hence $\psi(\alpha)\neq 1$ and therefore Theorem~\ref{thm:nontriviality} can be applied. Also note that since $\tilde{\pi}^{(n)}_r/\tilde{\pi}^{(n+1)}_r$ is $\mathbb{Z}$-torsion free, $\psi(\alpha^m)=1$ only if $m=0$. We claim that this implies that the kernel of
$$
\ov\phi:G\to \tilde{\pi}/\tilde{\pi}_r^{(n+2)}
$$
is contained in $G^{(1)}$. For suppose that $\ov\phi(\mu_K^mc)=1$ where $c\in G^{(1)}$. Then, since $i_*(G^{(1)})\subset \pi^{(n+1)}$, $G^{(1)}$ is contained in the kernel of
$$
\ov\psi:G\overset{\ov{\phi}}\to \tilde{\pi}/\tilde{\pi}_r^{(n+2)}\to \tilde{\pi}/\tilde{\pi}_r^{(n+1)},
$$
implying that $1=\ov{\psi}(c\mu_K^m)=\ov{\psi}(\mu_K)^m$. But since $\ov{\psi}(\mu_K)=\psi(\alpha)$, this is a contradiction unless $m=0$. Thus the kernel of $\ov{\phi}$ is contained in $G^{(1)}$ as claimed.

Now, by Theorem~\ref{thm:nontriviality}, if $P$ denotes the kernel of the map
$$
\mathcal{A}_0(K)\overset{i_*}{\to} H_1(M;\mathbb{Q}\Lambda)\overset{j_*}\to H_1(V;\mathbb{Q}\Lambda).
$$
then $P\subset P^\perp$ with respect to the classical Blanchfield form of $K$. Examine the commutative diagram below where $P$ is the kernel of the bottom horizontal composition. To justify the isomorphism in the bottom row, recall that $H_1(V;\mathbb{Q}\Lambda)$ is identifiable as the ordinary rational homology of the covering space of $V$ whose fundamental group is the kernel of $\psi:\tilde{\pi}\to \Lambda$. Since this kernel is precisely $\tilde{\pi}^{(n+1)}_r$, we have that
$$
H_1(V;\mathbb{Q}\Lambda)\cong (\tilde{\pi}^{(n+1)}_r/[\tilde{\pi}^{(n+1)}_r,\tilde{\pi}^{(n+1)}_r])\otimes_{\mathbb{Z}} \mathbb{Q}
$$
as indicated in the diagram
\begin{equation*}
\begin{CD}
G^{(1)}      @>i_*>>    \pi_1(M)^{(n+1)}  @>j_*>>   \tilde{\pi}^{(n+1)}_r  @>>>
\frac{\tilde{\pi}^{(n+1)}_r}{\tilde{\pi}^{(n+2)}_r} \\
  @VV{\pi}V   @VVV        @VVV       @VVjV\\
\mathcal{A}_0(K)     @>i_*>>  H_1(M;\mathbb{Q}\Lambda)    @>j_*>> H_1(V;\mathbb{Q}\Lambda) @>\cong>>
 \frac{\tilde{\pi}^{(n+1)}_r}{[\tilde{\pi}^{(n+1)}_r,\tilde{\pi}^{(n+1)}_r]}\otimes_{\mathbb{Z}} \mathbb{Q}.\\
\end{CD}
\end{equation*}
Since, by definition,
$$
\tilde{\pi}^{(n+2)}_r\equiv \text{kernel}(\tilde{\pi}^{(n+1)}_r\to (\tilde{\pi}^{(n+1)}_r/[\tilde{\pi}^{(n+1)}_r,\tilde{\pi}^{(n+1)}_r])\otimes_{\mathbb{Z}} \mathbb{Q}))
$$
the far-right vertical map $j$ is injective. Thus the kernel of the top horizontal composition is precisely $\pi^{-1}(P)$, which is precisely $\tilde{P}$. This identifies the image of the map $\ov{\phi}:G\to \pi/\pi^{(n+2)}_r$ as $G/\tilde{P}$ for a submodule $P\subset \mathcal{A}_0(K)$ where $P\subset P^\perp$. Thus $\rho(M_K,\ov\phi)$ is a first-order signature.

\noindent This completes the proof of Theorem~\ref{thm:Bingdouble}.
\end{proof}

In examining the proof above, one sees that we made little use of the fact that $T$ was a trivial link. Indeed, with slight modifications, the proof really establishes the following more general result. The more general result says that if one infects a slice link by a knot whose first-order signatures are large then the resulting link is not a slice link. This generalizes Harvey's ~\cite[Theorem 5.4]{Ha2} where it was shown under identical hypotheses that $\rho_0(K)$ obstructs $T(\alpha,K)$ from being slice.

\begin{thm}\label{thm:basiclink2} Suppose $T$ is a slice link of $m$ components, $n\geq 1$ and $\alpha$ is an unknotted circle in $S^3-T$ with $[\alpha]\in\pi_1(S^3-T)^{(n)}$ and $[\alpha]\notin\pi_1(M_T)^{(n+1)}_H$. Let $L$ denote $T(\alpha,K)$, the result of infection of $T$ along $\alpha$ using the knot $K$. If $L$ is topologically slice in a rational homology $4$-ball (or is even a rationally $(n+2)$-solvable link) then one of the first order signatures of $K$ is less in absolute value than the Cheeger-Gromov constant of $M_T$.
\end{thm}

\begin{proof}
%\begin{proof}[Proof of Theorem~\ref{thm:basiclink2}] 
The following modifications are necessary to the previous proof. We use the fact that $V$ is a (putative) rational $(n+2)$-solution to apply Theorem~\ref{thm:rho=0} when needed. Then instead of concluding that $\rho(M_K,\ov\phi)=0$ we have only that
$$
|\rho(M_K,\ov\phi)|=|\rho(M_T,\phi_T)|<C_{M_T}.
$$
\end{proof}

Before moving on to more general results, we give another application.

\vspace{.3in}

In ~\cite[Section 6]{Ha2} Harvey considered a filtration $\mathcal{F}^m_{
(n)}$ of the $m$-component string link concordance group wherein a string link $L$ is $(n)$-solvable if its closure $\hat{L}$ is an $(n)$-solvable link in the sense of ~\cite[Section 8]{COT}. The restriction of this filtration to boundary string links, $\mathcal{B}(m)$ was denoted $\mathcal{BF}^m_{(n)}$. Harvey defined specific real-valued higher-order signature invariants, $\rho_n$ of string links. She showed that \emph{each} $\rho_n$ gives a homomorphism $\boldsymbol{\rho_n}:\mathcal{B}(m)\to \mathbb{R}$, and induces a homomorphism
$$
\rho_n:\mathcal{BF}^m_{(n)}/\mathcal{BF}^m_{(n+1)}\to \mathbb{R}
$$
whose image, for any $m>1$, contains an infinite dimensional  vector subspace (over $\mathbb{Q}$) of $\mathbb{R}$. This was slightly improved to $\mathcal{BF}^m_{(n)}/\mathcal{BF}^m_{(n.5)}$ in ~\cite[Theorem 4.5]{CH2}. From this she concluded that (we incorporate the improvement of ~\cite[Theorem 4.5]{CH2})

\begin{thm}[{\cite[Theorem 6.8]{Ha2}}]\label{thm:shellythm} For any $m>1$ the abelianization of $\mathcal{BF}^m_{(n)}/\mathcal{BF}^m_{(n.5)}$ has infinite rank, and so $\mathcal{BF}^m_{(n)}/\mathcal{BF}^m_{(n.5)}$ is an infinitely generated subgroup of $\mathcal{F}^m_{(n)}/\mathcal{F}^m_{(n.5)}$.
\end{thm}

Our examples cannot be detected by any of Harvey's $\{\rho_n\}$ and so we can show that

\begin{cor}\label{cor:kernelrhon} For any $m>1$, $n\geq 2$, the kernel of Harvey's
$$
\rho_n:\mathcal{BF}^m_{(n)}/\mathcal{BF}^m_{(n.5)}\to \mathbb{R}
$$
contains an infinitely generated subgroup.
\end{cor}

\begin{proof}
%\begin{proof}[Proof of Corollary~\ref{cor:kernelrhon}] 
Let $\{K_i\}$ be an infinite set of Arf invariant zero knots such that $\{\rho_0(K_i)\}$ is a $\mathbb{Q}$-linearly independent subset of $\mathbb{R}$ (the existence of such a set was established in ~\cite[Proposition 2.6]{COT2}). Let $R_1$ be the ribbon knot $9_{46}$. It is easy to see, by taking a subset if necessary, that we can assume that $\{\rho_0(K_i), \rho^1(M_{R_1})\}$ is linearly independent. Let $J_i$ denote the knot of Figure~\ref{fig:examplehighersigs} with $K_1=K_2=K_i$. By ~\cite[Proposition 3.1]{COT2} $J_i$ is a $(1)$-solvable knot. Fix $m>1$ and let $T$ denote the trivial $m$-string link in $D^2\times I$. Fixing $n\geq 2$, choose a circle $\alpha\in F^{(n-1)}-F^{(n)}$, where $F$ is the group of the exterior of $T$, such that $\alpha$ bounds a disk in $D^2\times I$. Let $L_i$ denote $T(\alpha,J_i)$, the string link obtained by infecting $T$ along $\alpha$ using the knot $J_i$. The closure $\hat{L}_i$ is obtained from the trivial link  (which is $(n)$-solvable) by a $(1)$-solvable knot along a circle in $F^{(n-1)}$. Thus by Lemma~\ref{lem:nsolv}, $\hat{L}_i$ is $(n)$-solvable in the sense of ~\cite{COT}. Consequently $L_i\in \mathcal{F}^m_{(n)}$. It is easily seen that $L_i$ is a boundary string link (see ~\cite[Section 2]{CO2}), so
$$
L_i\in \mathcal{BF}^m_{(n)}.
$$
It follows directly from Harvey's formula ~\cite[Theorem 5.4]{Ha2} that $\rho_n(L_i)=0$ (indeed all of her $\rho_j$ vanish for these links). Consider the subgroup of $\mathcal{BF}^m_{(n)}$ generated by $\{L_i\}$. Suppose this were finitely generated. Then there is a subset $\{L_1,...,L_k\}$ that is a generating set. Consider $L_N$ for some $N>k$. Then the closure of the product $L=L_NL_{i_1}^{\epsilon_1}L_{i_2}^{\epsilon_2}...L_q^{\epsilon_q}$ is $(n.5)$-solvable for $i_j\in \{1,...,k\}$ and $\epsilon_j\in\{\pm 1\}$. A crucial point is now the observation that $\hat{L}$ can be obtained from the trivial link by multiple infections on curves $\alpha$ and $\alpha_i$, all lying in $F^{(n)}-F^{(n+1)}$, where the infection along $\alpha_N$ is done using $J_N$ and the other infections are done using copies of $J_1,...,J_k$ or their mirror images (if $\epsilon_j=-1$). The proof of Theorem~\ref{thm:basiclink} applies verbatim to this situation (although it was stated above for only one infection) because the crucial Theorem~\ref{thm:nontriviality} applies to the Alexander module of each infection knot separately. The conclusion is that some first first-order signature of $J_N$ is equal to some linear combination of first-order signatures of the knots $\{J_1,...,J_k\}$. We saw in Example~\ref{ex:first-ordersigs} that a first order signature for $J_i$ is an element of the set $\{\rho_0(K_i),\rho^1(R_1)+2\rho_0(K_i)\}$. It follows that $\rho_0(K_N)$ is a (possibly trivial) linear combination of $\{\rho_0(K_1),...,\rho_0(K_k),\rho^1(M_{R_1})\}$, contradicting our choice of $\{K_i\}$. Therefore the subgroup of $\mathcal{BF}^m_{(n)}$ generated by $\{L_i\}$ is infinitely generated.
\end{proof}

\section{Iterated Bing doubles and higher-order $L^{(2)}$-signatures}\label{sec:BingdoublesJn}

The techniques of the proof of Theorem~\ref{thm:Bingdouble} and Theorem~\ref{thm:basiclink2} can be generalized to include the iterated Bing doubles of more and more subtle knots, in particular knots whose classical signatures \emph{and} first-order signatures (and Casson-Gordon invariants) vanish. For specificity first consider the family of knots $J_n$ from Figure~\ref{fig:family}. If $n>1$ these have vanishing classical signatures, first-order signatures and vanishing Casson-Gordon invariants. Yet we show that higher-order signatures obstruct their iterated Bing doubles from being slice. For the family $J_n$, these higher-order signatures can be calculated, ``up to a constant,'' in terms of the classical signatures of $J_0$, so we formulate our results terms of $\rho_0(J_0)$ rather than in terms of an $n^{th}$ order signature of $J_n$. Since the proof will emphasize the structure of $J_n$ as obtained from $J_0$ by applying an $n$-fold doubling operator, we will use the notation $J_0=K$ and $J_n=J_n(K)$.

\begin{figure}[htbp]
\setlength{\unitlength}{1pt}
\begin{center}
\begin{picture}(230,117)
\put(0,0){\includegraphics{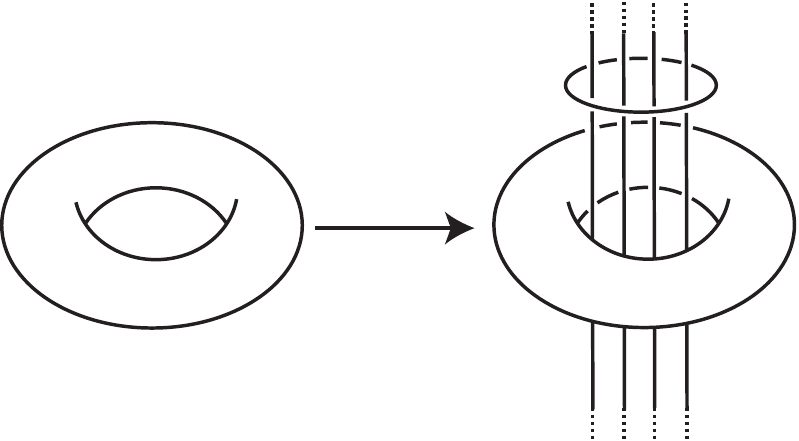}}
\put(150,100){$\tilde{\alpha}$}
\put(105,68){$g^{n-k}$}
\put(18,19){$S^3-J_{n-k}(K)$}
\put(210,18){$T$}
\end{picture}
\end{center}
\caption{$L_n(K)$}\label{fig:linkfamily}
\end{figure}

\begin{thm}\label{thm:mainlink} Suppose $T$ is a trivial link of $m$ components, $k$ and $n$ are positive integers such that $1\leq k\leq n$ and $\alpha$ is an unknotted circle in $S^3-T$ that represents an element in $F^{(k)}-F^{(k+1)}$ where $F=\pi_1(S^3-T)$, $K$ is a knot with Arf($K)=0$, and  $L_n(K)$ denotes $T(\alpha,J_{n-k}(K))$, the result of infection of $T$ along $\alpha$ using the knot $J_{n-k}(K)$ shown in Figure~\ref{fig:linkfamily}. Then there is a positive constant $C$ (independent of $K$) such that if $|\rho_0(K)|>C$, then $L_n(K)$ is not topologically slice in a rational homology ball. Moreover, $L_n(K)$ is $(n)$-solvable but not rationally $(n+1)$-solvable. If $L_n(K)$ is expressed as the closure of the $m$-component string link $\mathcal{SL}$ then no non-zero multiple of $\mathcal{SL}$ has closure that is rationally $(n+1)$-solvable.
\end{thm}

\begin{cor}\label{cor:BingdoubleJnK} For any $n$, there is a constant $C$ such that for any knot $K$ with Arf$(K)=0$ and $|\rho_0(K)|>C$ the Bing double of $J_{n-1}(K)$ is $(n)$-solvable but not slice nor even rationally $(n+1)$-solvable.
\end{cor}

\begin{proof}[Proof of Corollary~\ref{cor:BingdoubleJnK}] As we have seen in Figure~\ref{fig:bingeta},  a Bing double is obtained by a single infection of the trivial link of two components along a circle $\alpha$ representing the generator of the non-zero group $F^{(1)}/F^{(2)}$ where $F$ is the free group on two letters. The result then follows directly from Theorem~\ref{thm:mainlink} with $k=1$.
\end{proof}

\begin{cor}\label{cor:iteratedbing}Suppose $k$ and $n$ are positive integers. Then there is a constant $C$ such that if $K$ is any knot with Arf($K)=0$ such that $|\rho_0(K)|>C$, then the $k$-fold iterated Bing double of $J_{n-k}(K))$ is $(n)$-solvable but not slice nor even rationally $(n+1)$-solvable.
\end{cor}

\begin{proof}[Proof of Corollary~\ref{cor:iteratedbing}] As discussed in the proof of Theorem~\ref{thm:Bingdouble}, the $k$-fold iterated Bing double can be obtained from the trivial $2^k$ component link $T$ by a single infection, using the knot $J_{n-k}(K)$, along a circle $\alpha$ representing, in $\pi_1(M_T)$, an element of $F^{(k)}-F^{(k+1)}$. The result then follows directly from Theorem~\ref{thm:mainlink}.
\end{proof}

\begin{remark}[Remarks on Theorem~\ref{thm:mainlink}]

\

\begin{enumerate}

\item The restriction to Arf($K$)$=0$ is only to guarantee that $L_n(K)$ is $(n)$-solvable. It is not necessary for the conclusion that $L_n(K)$ is \emph{not} $(n+1)$-solvable.
\item  Using the different techniques of ~\cite[Theorem 9.1]{CHL3} one can show that $L_n(K)$ is not even rationally $(n.5)$-solvable, and one can choose $C$ independently of $n$ and $k$ (in fact $C$ can be chosen to be the Cheeger-Gromov constant of the zero surgery on the $9_{46}$ knot). We sketch the proof. Suppose that $L_n(K)$ were rationally $(n.5)$-solvable via $V$. Since $L=L_n(K)$ is obtained from the trivial link $T$ by infection along $\alpha$ using the knot $J_{n-k}(K)$, there is a cobordism $E$, as in Figure~\ref{fig:mickey}, with boundary components $M_L$, $M_T$ and $M_{J_{n-k}(K)}$. Cap off the $M_L$ boundary component using $V$ and cap off the $M_T$ boundary component using $\sharp_b S^1\times B^3$. The result, $W_0$, has boundary $M_{J_{n-k}(K)}$. Again by Figure~\ref{fig:mickey}, there is a cobordism $E'$ whose boundary components are $M_{J_{n-k}(K)}$, $M_R$ and two copies of $M_{J_{n-k-1}(K)}$. Adjoining $E'$ to $W_0$ we obtain $W_1$ whose boundary is $M_R$ and two copies of $M_{J_{n-k-1}(K)}$. Continuing in this way, we end up with a $4$-manifold, $W$, whose boundary is $n-k$ copies of $M_R$ and $2^{n-k}$ copies of $M_K$. With respect to a coefficient system $\pi_1(W)\to \pi_1(W)/\pi_1(W)^{(n+1)}_r$ the signature defects of the pieces of this cobordism are all zero, since these pieces are an $(n.5)$-solution, a slice disk exterior and many copies of the cobordisms $E$ of Section~\ref{signatures}. The signatures of the first two types vanish by ~\cite[Theorem 4.2]{COT} (see our Theorem~\ref{thm:rho=0}) and the signatures of $E$ vanish by Lemma~\ref{lem:mickeysig}. Therefore the sum of the $\rho$-invariants of the boundary components is zero. The sum of the contributions from the $M_R$ boundary components is bounded by a multiple of $C$. The $M_K$ components contribute some multiple of $\rho_0(K)$. A more careful analysis shows that in fact these multiples are comparable and one can conclude that $|\rho_0(K)|<C$. However this analysis depends crucially on a version of our Lemma~\ref{lem:selfannihil} under vastly weaker hypotheses. This is proved in ~\cite{CHL3}.

 \item If we use a different family $J_n(K)$ as shown in Figure~\ref{fig:betterfamilyJ_n} ($T$ is the trefoil knot) patterned after the ribbon knot $R$, obtained by setting $J_{n-1}=U$ in Figure~\ref{fig:betterfamilyJ_n}, then much better results are possible. The key difference is that $\rho^1(R)\neq 0$ by an analysis as in Example~\ref{ex:nonzerofirst-ordersigs}. In particular, applying the techniques of ~\cite[Theorem 9.1]{CHL3} to $L_n(K)$ (defined using this different family) one can completely eliminate the Cheeger-Gromov constant (replace it by $C=0$). In terms of the proof sketch just above it allows us to cap off all of the copies of $M_R$. In this way we get \emph{specific} examples of $(n)$-solvable knots none of whose iterated Bing doubles is slice.

\begin{figure}[htbp]
\begin{center}
\begin{picture}(138,135)
\put(0,0){\includegraphics{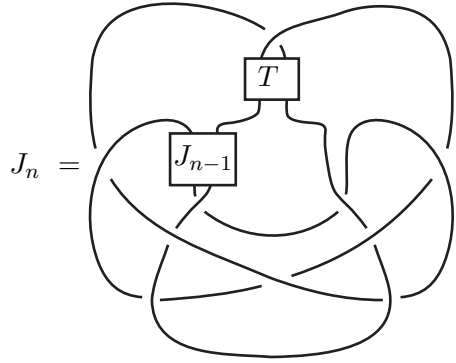}}
%\put(79,136){}
\put(64,103){$T$}
\put(32,73){$J_{n-1}$}
\put(-30,70){$J_{n}~=$}
\end{picture}
\end{center}
\caption{A different family of $(n)$-solvable knots $J_n$}\label{fig:betterfamilyJ_n}
\end{figure}

\end{enumerate}
\end{remark}

\begin{proof}[Proof of Theorem~\ref{thm:mainlink}] The proof is not substantially different from  that of Theorem~\ref{thm:basiclink}, but is notationally much more complicated. Without loss of generality we can assume that $L\equiv L_n(K)$ is the closure of a string link $\mathcal{SL}$ as in the last clause of the theorem. The closure of a multiple of $\mathcal{SL}$ is just a particular connected-sum of copies of $L$. Hence it suffices to show that, if $|\rho_0(K)|$ is sufficiently large, then $\#^M_{j=1}L$ in \emph{not} rationally $(n+1)$-solvable for any $M>0$.

We now state one lemma and two theorems. Assuming these three, we finish the proof of Theorem~\ref{thm:mainlink}. The rest of the section will then be devoted to the proofs of these three results.

We first claim that $L$ can be obtained from a ribbon link by multiple infections along curves lying in the $n^{th}$-derived subgroup of the ribbon group. Specifically let $U$ be the unknot, let $R_i\equiv J_i(U)$ denote the family of ribbon knots obtained recursively by setting $J_0=U$ in Figure~\ref{fig:family} or by applying the $9_{46}$ operator $n$ times to $K=U$ as in Figure~\ref{fig:Rdoubling}. Then $L_n(U)=T(\alpha,J_{n-k}(U))=T(\alpha,R_{n-k})$ is a ribbon link. The precise definition of the circles
$\alpha^{n-k}_*$ (clones) will be given in the proof of Lemma~\ref{cor:linknsolvable}.

\begin{lem}\label{cor:linknsolvable} $L_n(K)$ can be obtained from the slice link $L_n(U)=T(\alpha,R_{n-k})$ as the result of $2^{n-k}$ infections using the knot $K$ each time, along the clones $\alpha^{n-k}_*$ that lie in $\pi_1(S^3-L_n(U))^{(n)}$.
\end{lem}

\begin{thm}\label{thm:iteratedlink} Let $T(\alpha,R_{n-k})$ be as above. Suppose $W$ is an \emph{arbitrary} rational $(n)$-solution for $M_{T(\alpha,R_{n-k})}$. Then for at least one of the $2^{n-k}$ clones $\alpha^{n-k}_*$, $j_*(\alpha^{n-k})\notin \pi_1(W)_r^{(n+1)}$ where
$$
j_*:\pi_1(M_{T(\alpha,R_{n-k})})\to \pi_1(W)
$$
is induced by inclusion.
\end{thm}

\begin{thm}\label{thm:linkinfection} Let $R$ be a slice link of $m$ components
($n\ge1$) and $M_R$ the $0$-framed surgery on $R$. Suppose there
exists a collection of homotopy classes
\[
[\eta_i]\in\pi_1(M_R)^{(n)}, \quad 1\leq i\leq N,
\]
that has the following property: For \emph{any}
rational $(n)$-solution $W$ of $M_R$ there exists \emph{some} $i$ such that
$j_*(\eta_i)\notin\pi_1(W)^{(n+1)}_r$ where $j_*:\pi_1(M_R)\to\pi_1(W)$.

Then, for any oriented trivial link $\{\eta_1,\dots,\eta_N\}$ in $S^3\smallsetminus R$
that represents the $[\eta_i]$, and for any $N$-tuple
$\{K_1,\dots,K_N\}$ of Arf invariant zero knots for which
either each $\rho_0(K_i)>C_{M_R}$ (the Cheeger-Gromov constant of $M_R$), or each $\rho_0(K_i)<-C_{M_R}$, the link
\[
L=R(\eta_1,...,\eta_N,K_1,...,K_N)
\]
is $(n)$-solvable but no positive multiple of it is slice (nor even rationally $(n+1)$-solvable). (If the Arf invariant condition is dropped then $L$ is merely rationally $n$-solvable).
\end{thm}

Now, assuming, Lemma~\ref{cor:linknsolvable} and Theorems~\ref{thm:iteratedlink} and ~\ref{thm:linkinfection}, we finish the proof of Theorem~\ref{thm:mainlink}. We claim that Theorem~\ref{thm:linkinfection} applies to $L=L_n(K)$, and that Lemma~\ref{cor:linknsolvable} and Theorem~\ref{thm:iteratedlink} merely serve to show that $L$ satisfies the hypotheses of Theorem~\ref{thm:linkinfection}. Specifically, we seek to apply Theorem~\ref{thm:linkinfection} with $R=T(\alpha,R_{n-k})$, $N=2^{n-k}$, $K_i=K$ for all $i$, $L=L_n(K)$, and $\{\eta_1,\dots,\eta_N\}=\{\alpha^{n-k}_*\}$. To verify the hypotheses of Theorem~\ref{thm:linkinfection} we need that: $R$ is a slice link, and that each $\alpha^{n-k}_*\in \pi_1(M_R)^{(n)}$, both of which are guaranteed by Lemma~\ref{cor:linknsolvable}. Moreover we need that
for any rational $(n)$-solution $W$ for $M_R$ there exists \emph{some} clone $\alpha^{n-k}$ such that
$j_*(\alpha^{n-k})\notin\pi_1(W)^{(n+1)}_r$. But this is guaranteed by Theorem~\ref{thm:iteratedlink}. Therefore we have verified the hypotheses of Theorem~\ref{thm:linkinfection}.  By the conclusion of that theorem, $L$ is $(n)$-solvable but no positive multiple of it is slice (nor even rationally $(n+1)$-solvable), as long as $|\rho_0(K)|>C$ where $C$ is the Cheeger-Gromov constant for $M_{T(\alpha,R_{n-k})}$. This concludes the proof of Theorem~\ref{thm:mainlink}.
\end{proof}

Now we turn to the proof of Theorem~\ref{thm:linkinfection}, which is a very general analog, for links, of ~\cite[Theorem 4.2]{CT} (for knots).

\begin{proof}[Proof of Theorem~\ref{thm:linkinfection}] Supposing that such $R$ and $\eta_i$ exist, let $L=R(\eta_1,...,\eta_N,K_1,...,K_N)$ for knots $K_i$ such that, for each $i$, Arf($K_i)=0$ and $\rho_0(K_i)>C_{M_R}$ (the Cheeger-Gromov constant of $M_R$).

Since $L$ is, by hypothesis, the result of infections on an $(n)$-solvable link (in fact a slice link) along circles lying in the $n^{th}$ term of the derived series $L$ is $(n)$-solvable (merely rationally $(n)$-solvable without the Arf invariant condition) by ~\cite[Proposition 3.1]{COT2}.

For the remainder of the proof, we proceed by contradiction. Suppose that $\tilde L\equiv\#^p_{j=1}L$ were rationally $(n+1)$-solvable for some $p>0$. Then there would exist a rational $(n+1)$-solution $V$ with $\partial V= M_{\tilde L}$, the zero framed surgery on $\tilde L$. Using this we construct a particular rational $(n)$-solution $W$ for $M_{R}$ as follows (shown schematically in Figure~\ref{fig:cobordismlink}).
\begin{figure}[htbp]
\setlength{\unitlength}{1pt}
\begin{center}
\begin{picture}(326,128)
\put(155, 43){$C$}
\put(39, 111){$Z^1_1$}
\put(71, 111){$Z^1_N$}
\put(39, 80){$E^1$}
\put(160, 80){$E^2$}
\put(278, 80){$E^p$}
\put(161, 111){$Z^2_1$}
\put(193, 111){$Z^2_N$}
\put(-20,62){$M_{L}^1$}
\put(333,62){$M_{L}^p$}
\put(-20,97){$M_{R}$}
\put(105,97){$M_{R}^2$}
\put(223,97){$M_{R}^p$}
\put(126, 108){$Y^2$}
\put(245, 108){$Y^p$}
\put(278, 111){$Z^p_1$}
\put(310, 111){$Z^p_N$}
\put(155,9){$V$}
\put(-20,23){$M_{\tilde{L}}$}
\put(0,0){\includegraphics{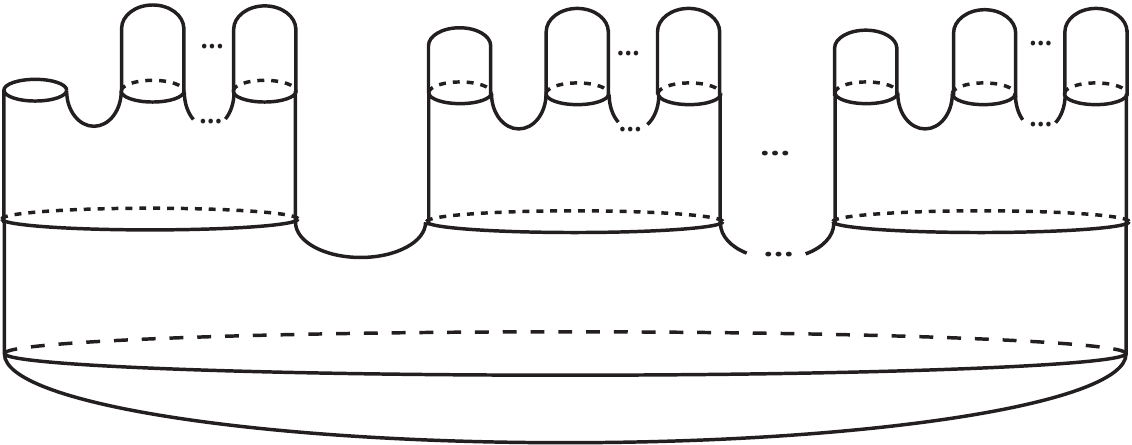}}
\end{picture}
\end{center}
\caption{The rational $(n)$-solution $W$ for $M_{R}$} \label{fig:cobordismlink}
\end{figure}
Here $C$ is the standard cobordism from $M_{\tilde L}$ to the disjoint union of $p$ copies of $M_{L}$. This cobordism is discussed in detail in ~\cite[Section 4]{COT2}. Cap off the boundary component $M_{\tilde L}$ using the rational $(n+1)$-solution $V$. Since $L$ is obtained from the link $R$ by infection on circles $\eta_i$ using the knots $K_i$, there is a cobordism $E$, as shown in Figure~\ref{fig:mickey}, such that
$$
\partial E= -M_L\sqcup M_R \sqcup_{i=1}^{N} M_{i}
$$
where we abbreviate $M_{K_i}$ by $M_i$. Add a copy of $E$ to each of the $p$ copies of $M_{L}$. We denote these copies by $E^j, 1\leq j\leq p$. Now, for each $i$, cap off each of the $p$ copies of $M_{i}$ with a $(0)$-solution $Z^j_i$ for $K_i$ (we can assume that $\pi_1(Z^j_i)=\mathbb{Z}$ by ~\cite[p.108]{COT2}~\cite[Appendix 5]{COT2}) and cap off each of the copies of $M_R$, except the ``first,'' with a copy, $Y^j, 2\leq j\leq p$, of the exterior $Y$ of a set of slicing disks for the slice link $R$. The resulting manifold $W$ then has a single copy of $M_{R}$ as its boundary.

\begin{lem}\label{lem:H2} $W$ is a rational $(n)$-solution for $M_R$.
\end{lem}

\begin{proof}[Proof of Lemma~\ref{lem:H2}] By Definition~\ref{defn:rationalnsolvable}, we must show that
\begin{itemize}
\item [(1)] $H_1(M_R;\mathbb{Q})\to H_1(W;\mathbb{Q})$ is an isomorphism, and
\item  [(2)] $W$ admits a rational $(n)$-Lagrangian with rational $(n)$-duals.
\end{itemize}
First we claim that
$$
H_2(W;\mathbb{Q})\cong H_2(V;\mathbb{Q})\oplus_{i,j} H_2(Z_i^j;\mathbb{Q}).
$$
Since $V$ is a rational $(n+1)$-solution for $M_{\tilde L}$, the inclusion-induced map
$$
j_*:~H_1(M_{\tilde L};\mathbb{Q})\to H_1(V;\mathbb{Q})
$$
is an isomorphism. It follows from duality that
$$
j_*:~H_2(M_{\tilde L};\mathbb{Q})\to H_2(V;\mathbb{Q})
$$
is the zero map. Therefore if we examine the Mayer-Vietoris sequence with $\mathbb{Q}$-coefficients,
$$
H_2(C)\oplus H_2(V)\overset{\pi_*}{\lra} H_2(C\cup V)\to H_1(M_{\tilde L})\overset{(i_*,j_*)}{\longrightarrow}H_1(C)\oplus H_1(V),
$$
we see that $\pi_*$ induces an isomorphism
$$
(H_2(C)/(i_*(H_2(M_{\tilde L})))\oplus H_2(V)\cong H_2(V\cup C).
$$
Moreover recall that $C$ is obtained from a collar of the disjoint union of $p$ copies of $M_{L}$ by adding $(p-1)$ $1$-handles (to connect the components) and then adding $m(p-1)$ $2$-handles that have the effect of equating pairwise the meridional elements of the copies $L$. In this way we see that, for any of the boundary components $M_L$, $H_1(M_L;\mathbb{Q})\cong H_1(C;\mathbb{Q})\cong \mathbb{Q}^m$ generated by a set of meridians, and that $H_2(C;\mathbb{Q})\cong\oplus_{j=1}^p H_2(M_L;\mathbb{Q})$ (this is analyzed in more detail in ~\cite[p. 113-114]{COT2}). It is easy to see that a basis of $i_*(H_2(M_{\tilde L}))$ is formed from the sum, $1\leq j\leq p$ of the elements of natural bases for each $H_2(M_L;\mathbb{Q})$. Thus
$$
H_2(V\cup C;\mathbb{Q})\cong   H_2(V;\mathbb{Q})\oplus(\oplus_{j=1}^p H_2(M_L;\mathbb{Q}))/D
$$
where $D\cong \mathbb{Q}^m$ is the diagonal subgroup. Now, recall that we have analyzed the homology of $E$ in Lemma~\ref{lem:mickeyfacts} and found that,
$$
H_1(M_{L})\overset{i_*}{\longrightarrow} H_1(E)
$$
is an isomorphism. Therefore the following Mayer-Vietoris sequence with $\mathbb{Q}$-coefficients is exact,
$$
\oplus^p_{j=1} H_2(M_{L}^j)\to \oplus^p_{j=1} H_2(E^j)\oplus H_2(V\cup C)\overset{\pi_*}{\to} H_2(V\cup C\sqcup^p_{j=1} E^j)\to 0.
$$
Moreover, from property $(3)$ of Lemma~\ref{lem:mickeyfacts},
$$
H_2(E)\cong \oplus^N_{i=1}H_2(M_i)\oplus H_2(M_R)
$$
where the latter $H_2(M_R)\cong H_2(M_L)$ in $H_2(E)$. Combining these facts we have that
$$
H_2(V\cup C\sqcup^p_{j=1} E^j)\cong H_2(V)\oplus_{j=1}^p\oplus_{i=1}^N H_2(M_i^j)\oplus_{j=1}^p(H_2(M_R^j)/D).
$$
The next step in the formation of $W$ was the addition of the slice exteriors $Y^j$ to the copies $M_R^j$ for $2\leq j\leq p$. Since $H_1(\partial Y^j)\to H_1(Y^j)$ is an isomorphism and $H_2(Y^j)=0$, the effect on $H_2$ of adding the $Y^j$ is merely to kill all the $H_2$ carried by the boundaries $H_2(M^j_R)$, $2\leq j\leq p$. Taking into account the diagonal relation, we have
\begin{equation}\label{eq:bordism}
H_2(V\cup C\cup E^j\cup Y^j)\cong H_2(V)\oplus_{j=1}^p\oplus_{i=1}^N H_2(M_i^j).
\end{equation}
The final step in the formation of $W$ was the addition of the $(0)$-solutions $Z^j_i$ to all the copies $M_i^j$ of $M_{K_i}$. Since, $Z^j_i$ is a $(0)$-solution, $H_1(M_i^j)\to H_1(Z_i^j)$ is an isomorphism and by duality $H_2(M_i^j)\to H_2(Z_i^j)$ is the zero map. Thus the effect on $H_2$ of adding the $Z_i^j$ is merely to kill all the generators of the $H_2(M_i^j)$ summand and add $H_2(Z_i^j)$. Thus we have
$$
H_2(W;\mathbb{Q})\cong H_2(V;\mathbb{Q})\oplus_{i,j}H_2(Z_i^j)
$$
This establishes the claim. Combining some of the observations above it also follows that $H_1(M_R;\mathbb{Q})\to H_1(W;\mathbb{Q})$ is an isomorphism.

We return now to the proof that $W$ is a rational $(n)$-solution for $M_R$. Since $V$ is a rational $(n+1)$-solution, it is a rational $(n)$-solution. Let $\{\ell_1,\dots,\ell_g\}$ be a collection of $(n)$-surfaces generating a rational $(n)$-Lagrangian for $V$ and $\{d_1,\dots,d_g\}$ be a collection of $(n)$-surfaces generating its rational $(n)$-duals. By definition, $2g=\text{rank}_\mathbb{Q}H_2(V;\mathbb{Q})$. Similarly, for each $i$ and $j$ take a collection of such $(0)$-surfaces $\{l^{ij}_1,..,l^{ij}_k\}$, $\{d^{ij}_1,..,d^{ij}_k\}$ for the $(0)$-solutions $Z_i^j$. Now taking these surfaces for $V$ together with the collections of surfaces for the $Z_i^j$, these collections have the required \emph{cardinality} (by the first part of the lemma) to generate a rational $(n)$-Lagrangian with rational $(n)$-duals for $W$. Since $\pi_1(V)^{(n)}$ maps into $\pi_1(W)^{(n)}$, the $(n)$-surfaces for $V$ are also $(n)$-surfaces for $W$. We need to show that the $(0)$-surfaces for $Z_i^j$ are $(n)$-surfaces for $W$.

The group $\pi_1(Z_j^i)\cong \mathbb{Z}$ is generated by the meridian of the knot $K^j_i$ in $M^j_i$.  This meridian is isotopic in $E_j$ to the infection curve $\eta^j_i\in M_R^j$. By hypothesis,
\begin{equation}\label{eq:nth}
[\eta^j_i]\in \pi_1(M_R^j)^{(n)}.
\end{equation}
Therefore
$$
j_*(\pi_1(Z^i_j))\subset \pi_1(W)^{(n)}.
$$
Hence \emph{any} surface in $Z_i^j$ is an $(n)$-surface for $W$. Moreover, by functoriality of the intersection form with twisted coefficients these collections of surfaces have the required intersection properties to generate a rational $(n)$-Lagrangian with rational $(n)$-duals for $W$. Hence $W$ is a rational $(n)$-bordism for $M_R$, as was claimed.  This completes the proof of Lemma~\ref{lem:H2}.
\end{proof}

It also follows from ~(\ref{eq:bordism}) just above that:

\begin{cor}\label{cor:n+1bordism} Let $X=V\cup C\cup E^j\cup Y^j$ so that $\partial X= M_R\cup_{i,j} M_i^j$. Then
the cokernel of
$$
H_2(\partial X;\mathbb{Q})\to H_2(X;\mathbb{Q})
$$
is $H_2(V;\mathbb{Q})$.
\end{cor}

We continue with the proof of Theorem~\ref{thm:linkinfection}. Now set $\Gamma= \pi_1(W)/\pi_1(W)^{(n+1)}_r$. Let $\psi:\pi_1(W)\to\Gamma$ be canonical surjection. Let $\phi:\pi_1(M_{R})\to\Gamma$ be the composition $\psi\circ j_*$. Thus by the hypothesis of Theorem~\ref{thm:linkinfection} there exists \emph{some} $i$ such that
$\phi(\eta_i)\neq 1$. We shall now compute $|\rho(M_{R},\phi)|$ using $W$, and find it to be greater than $C_R$. This contradiction will show that in fact $\tilde L\equiv\#^p_{j=1}L$ is not rationally $(n+1)$-solvable.

By definition we have,
$$
\rho(M_{R},\phi)= \sigma^{(2)}_\Gamma(W,\psi)-\sigma(W).
$$
By the additivity of the non Neumann and the ordinary signatures (\cite[Lemma 5.9]{COT}) the latter signatures are the sums of the corresponding signatures for the submanifolds $X$ and $Z^j_i$.

First consider $X$. Using Corollary~\ref{cor:n+1bordism} and the fact that $V$ is a rational $(n+1)$-solution, $X$ is what is called a \emph{rational $(n+1)$-bordism} in ~\cite[Section 5]{CHL3}. A rational $(n+1)$-bordism is similar to a rational $(n+1)$-solution except that its boundary need not be connected and the inclusion-induced maps on $H_1$ from its boundary components are unrestricted. Since $\Gamma^{(n+1)}=1$, by ~\cite[Theorem 5.9]{CHL3},
$$
\sigma^{(2)}_\Gamma(X)-\sigma(X)=0,
$$
as long as each of the boundary components, $M$, of $X$ satisfies the following alternative: either the induced coefficient system is trivial on $\pi_1(M)$, or
\begin{equation}\label{eq:alter}
\text{rank}_{\mathcal{K}\G}H_1(M;\mathcal{K}\G)=\beta_1(M)-1.
\end{equation}
This alternative is always satisfied if $\beta_1(M)=1$ (by ~\cite[Proposition 2.11]{COT}), as is the case for each $M_i^j$. That leaves only $M_R$ to consider. Let $B=\pi_1(X)$. We claim that there is a basis of $H_2(X;\mathbb{Q})$ consisting of surfaces $\Sigma\to X$ for which $\pi_1(\Sigma)\subset B^{(n+1)}$, which is what we call a $B^{(n+1)}$-surface. Recall from ~(\ref{eq:bordism}) that $H_2(X;\mathbb{Q})$ is generated by $H_2(V;\mathbb{Q})$ and by the $H_2(M^j_i;\mathbb{Q})$. Since $V$ is a rational $(n+1)$-solution, $H_2(V;\mathbb{Q})$ is generated by $\pi_1(V)^{(n+1)}$-surfaces, which are, a fortiori, $B^{(n+1)}$-surfaces. $H_2(M^j_i;\mathbb{Q})$ is generated by a capped-off Seifert surface for the knot $K_i^j$. Any circle on this Seifert surface lies in $\pi_1(M^j_i)^{(1)}$ and hence lies in $B^{(n+1)}$ since the meridian of $M_i^j$ lies in $B^{(n)}$ as we saw in ~(\ref{eq:nth}). Thus the Seifert surface is also a $B^{(n+1)}$-surface. This completes the verification of the claim. Choose a free group $F$ and a map $F\to \pi_1(M_R)$ inducing an isomorphism on $H_1$. Now consider the maps
$$
F\overset{i}\to \pi_1(X)\to B\overset{\psi}\to \G.
$$
Note that each of these maps induces isomorphisms on $H_1(-;\mathbb{Q})$. Now ~\cite[Proposition 2.11]{CH2} applies to both $F\to B$ and $\pi_1(M_R)\to B$, since $H_2(X;\mathbb{Q})$ has a basis of ker($\psi$)-surfaces since $B^{(n+1)}\subset \text{ker}(\psi)$. Thus
$$
H_1(F;\mathcal{K}\G)\cong H_1(M_R;\mathcal{K}\G)\cong H_1(B;\mathcal{K}\G).
$$
The rank of the first of these three is known to be $\beta_1(M)-1$ \cite[Lemma 2.12]{COT}. This completes the verification that $M_R$ satisfies the alternative ~(\ref{eq:alter}) and hence completes the verification that the the $\G$-signature defect of $X$ vanishes.

Now consider the $Z^j_i$. Let $\psi^j_i$ denote the restriction of $\psi$ to $\pi_1(Z_i^j)$. Then, by definition
$$
\sigma^{(2)}_\Gamma(Z_i^j)-\sigma(Z_i^j)=\rho(M_i^j,\psi_i^j).
$$
However, since $\pi_1(Z^j_i)\cong \mathbb{Z}$, $\psi^j_i$  factors through $\mathbb{Z}$. Hence by properties $(2)$, $(3)$ and $(4)$ of Proposition~\ref{prop:rho invariants}
$$
\rho(M_i^j,\psi_i^j)=\rho_0(K_i)
$$
if $\psi^j_i(\eta^j_i)\neq 1$ and is zero if  $\psi^j_i(\eta^j_i)=1$. Note that here we have used the fact that the infection circle $\eta^j_i$ (in $M_R^j$) is isotopic (in $E_j$) to the meridian of $K_i^j$ in $M^j_i$ (see property $(4)$ of Lemma~\ref{lem:mickeyfacts}).

Putting all of these together we have
$$
\rho(M_{R},\phi)=\sum_{i=1}^N d_i\rho_0(K_i)
$$
where $d_i$ is the number of values of $j$ for which $\psi(\eta^j_{i})\neq 1$. Since our hypothesis is that either for each $i$
$$
\rho_0(K_i)> C_{M_R},
$$
or for each $i$
$$
\rho_0(K_i)< -C_{M_R},
$$
this is a contradiction unless each $d_i=0$. However recall $W$ is a rational $(n)$-solution for $M_R$ by Lemma~\ref{lem:H2}. Thus by hypothesis there exists \emph{some} $i$ such that
$j_*(\eta^1_i)\notin\pi_1(W)^{(n+1)}_r$ where $j_*:\pi_1(M_R)\to\pi_1(W)$. Hence for some $i$,
$$
\psi^j_i(\eta^1_i)\neq 1,
$$
so $d_i>0$. This is a contradiction, completing the proof of Theorem~\ref{thm:linkinfection}.
\end{proof}

Now we turn to the proof of Lemma~\ref{cor:linknsolvable}. To accomplish this we will show that $L_n(K)$ has a variety of different descriptions due to its ``fractal'' nature.

\begin{proof}[Proof of Lemma~\ref{cor:linknsolvable}] Recall $U$ denotes the trivial knot, and $J_0(K)\equiv K$. First we establish that $J_n(K)$ has an alternative description as the result of $2^n$ infections on the ribbon knot $R_n=J_n(U)$ using the knot $K$ as the infecting knot each time, along curves that lie in $\pi_1(S^3\backslash R_n)^{(n)}$. This will be established as part of a much more general result that says that $J_n(K)$ has many alternative descriptions.

To this end note that if $K$ is the trivial knot $U$ then it is easily seen by induction that each $J_n(U)$ is a ribbon knot that we denote $R_n$, $n\ge 0$, as shown in Figure~\ref{fig:ribbonfamily} (set $R_0=U$). For, if $R_{n-1}$ is a ribbon knot then $2$ parallels of it form a $2$-component ribbon link. Then $R_n$ is formed from this ribbon link by fusing together the $2$ components using a knotted band.

\begin{figure}[htbp]
\setlength{\unitlength}{1pt}
\begin{center}
\begin{picture}(143,151)
\put(0,0){\includegraphics{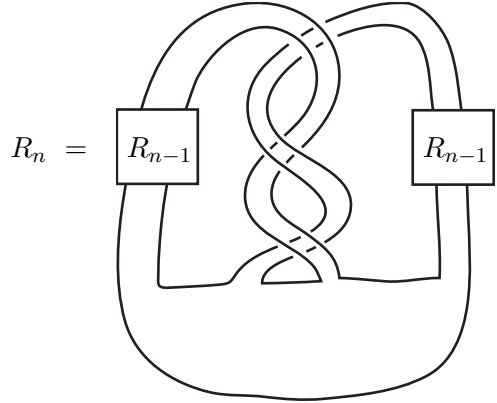}}
\put(-40,92){$R_{n}~=$}
\put(4,92){$R_{n-1}$}
\put(116,92){$R_{n-1}$}
\end{picture}
\end{center}
\caption{The recursive family of ribbon knots $R_{n}$}\label{fig:ribbonfamily}
\end{figure}

\noindent Now note that, for each $1\le i\le n$, because of the alternative description of infection as described in Section~\ref{sec:Introduction}, there are two inclusion maps
$$
f_\pm^{i}: S^3-R_{i-1}\ra S^3- R_{i}
$$
as suggested by Figure~\ref{fig:mapsf}.

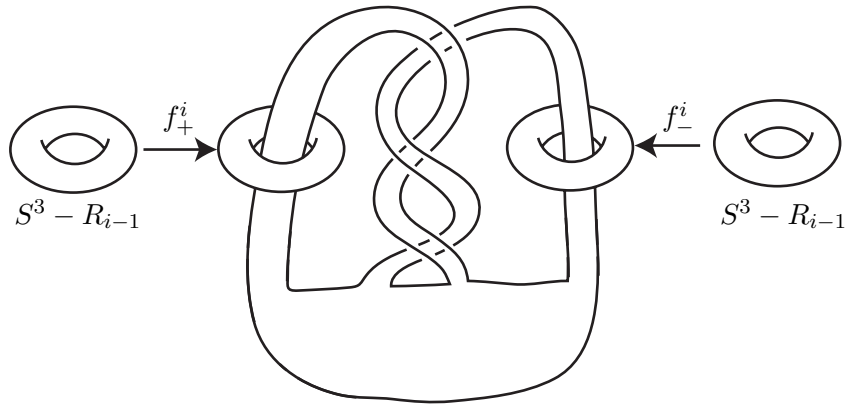
\begin{figure}[htbp]
\setlength{\unitlength}{1pt}
\begin{center}
\begin{picture}(314,150)
\put(0,0){\includegraphics{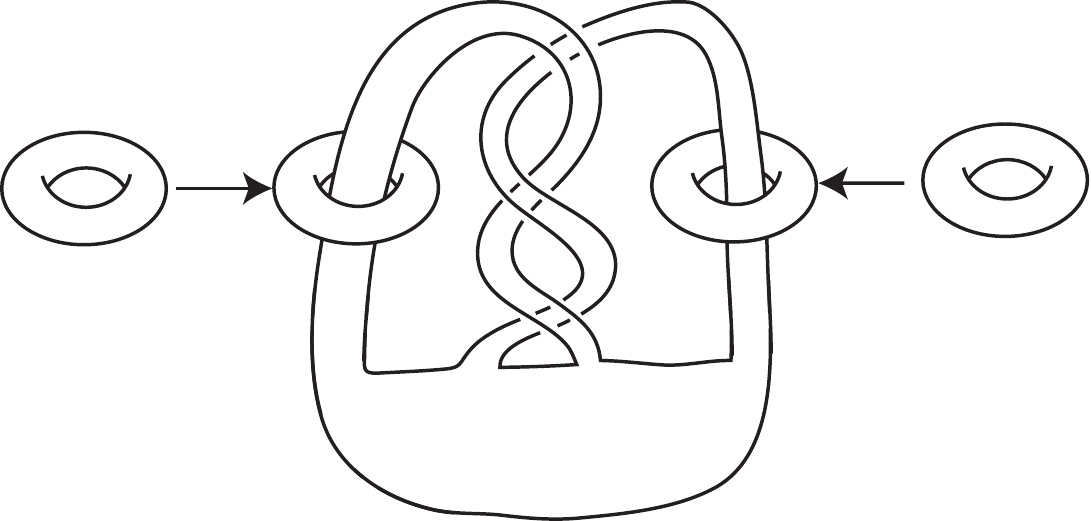}}
\put(58,105){$f^{i}_+$}
\put(247,105){$f^{i}_-$}
\put(2,68){$S^3-R_{i-1}$}
\put(269,68){$S^3-R_{i-1}$}
\end{picture}
\end{center}
\caption{The embeddings  $S^3-R_{i-1}\hookrightarrow S^3-R_{i}$}\label{fig:mapsf}
\end{figure}

Let $\eta^0$ denote the meridian of $R_0$, the trivial knot. Let $\eta^1_+, \eta^1_-$ denote the two images $f_\pm^1(\eta^0)$ in $S^3-R_1$. We call these \textbf{clones} of $\eta^0$. More generally, let $\{\eta^i_*\}$ denote the set of $2^i$ images of $\eta^0$ under the $2^i$ compositions $f_\pm^{i}\circ\dots\circ f_\pm^1$. Note that the induced maps
$$
(f_\pm^i)_*: \pi_1(S^3\backslash R_{i-1})\ra\pi_1(S^3\backslash R_{i})
$$
have images contained in the commutator subgroup. Thus the composition
$$
(f_\pm^{i})_*\circ\dots\circ(f_\pm^1)_*:  \pi_1(S^3\backslash R_0)\ra\pi_1(S^3\backslash R_1)^{(1)}\ra\dots\ra\pi_1(S^3\backslash R_i)^{(i)}
$$
has image in $\pi_1(S^3\backslash R_i)^{(i)}$. Therefore we see that each of the clones $\{\eta^i_*\}$ lies in $\pi_1(S^3\backslash R_i)^{(i)}$ and in particular each of the clones $\{\eta_*^n\}$ lies in $\pi_1(S^3\backslash R_n)^{(n)}$. The superscript $i$ of $\{\eta^i_*\}$ can serve to remind the reader in which term of the derived series it lies.

The following establishes that $J_n(K)$ has a variety of different descriptions.

\begin{prop}\label{prop:altdescriptions} For any knot $K$ and $i$, $0\leq i \leq n$, $J_n(K)$ can be obtained from $R_i$ by multiple infections along the $2^i$ clones
$$
\{\eta^i_*\}= ~\{f^{i}_{\pm}\circ\dots\circ f^1_{\pm}(\eta^0)\},
$$
using knot $J_{n-i}(K)$ as the infecting knot in each case, and each clone $\eta^i_*$ lies in $\pi_1(S^3-R_i)^{(i)}$.
\end{prop}

\begin{proof}[Proof of Proposition~\ref{prop:altdescriptions}] We proceed by induction on $i$. In the base case, $i=0$, for any $n$, there is only one clone, namely $\eta^0$ itself. Then the claim is merely that if one infects the unknot by $J_n(K)$ along a meridian then the result is $J_n(K)$, which is obviously true.

Assume that the proposition is true for some fixed $i-1$ for \emph{any} $n$ such that $n\geq i-1$. Then consider fixed $i$ and arbitrary $n$ subject to $n\geq i$. Recall that $S^3-J_n(K)$ can be obtained by deleting the two solid tori as shown in the Figure~\ref{fig:family2} and replacing them with two copies of $S^3-J_{n-1}(K)$.
\begin{figure}[htbp]
\setlength{\unitlength}{1pt}
\begin{center}
\begin{picture}(314,150)
\put(0,0){\includegraphics{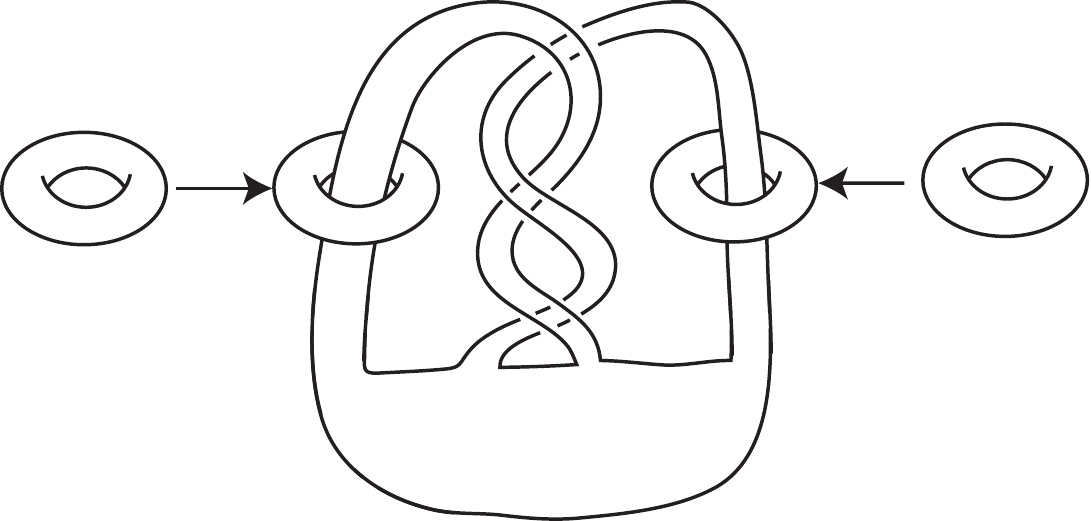}}
%\put(77,136){}
%\put(323,137){}
\put(2,68){$S^3-J_{n-1}$}
\put(269,68){$S^3-J_{n-1}$}
%\put(95,71){}
%\put(306,71){}
%\put(122,11){}
\end{picture}
\end{center}
\caption{One definition of $S^3-J_{n}$}\label{fig:family2}
\end{figure}

By the inductive hypothesis for $(n-1,i-1)$,  $S^3-J_{n-1}$ can be obtained from $S^3-R_{i-1}$ by infections on the $2^{i-1}$ clones $\{\eta^{i-1}_*\}\equiv ~\{f^{i-1}_{\pm}\circ\dots\circ f^1_{\pm}(\eta^0)\}$ (shown schematically by the very small solid tori in Figure~\ref{fig:alternativeviews} ) using the knot $J_{n-i}(K)$ as the infecting knot in each case. Thus replacing the $2^i$ solid tori shown in Figure~\ref{fig:alternativeviews} by copies of $S^3-J_{n-i}(K)$ yields $S^3-J_n$.
\begin{figure}[htbp]
\setlength{\unitlength}{1pt}
\begin{center}
\begin{picture}(314,150)
\put(0,0){\includegraphics{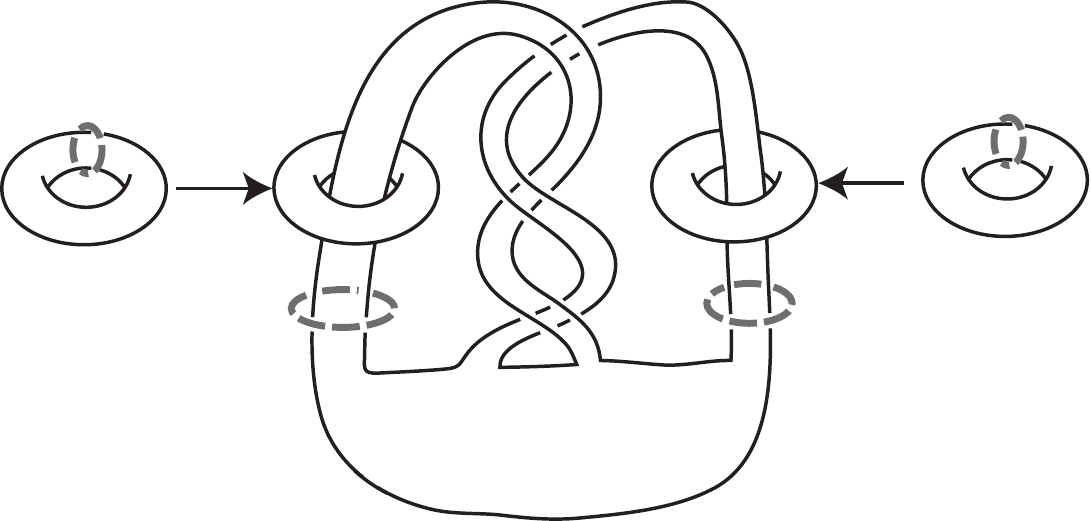}}
\put(58,105){$f^{i}_+$}
\put(247,105){$f^{i}_-$}
\put(2,68){$S^3-R_{i-1}$}
\put(269,68){$S^3-R_{i-1}$}
\end{picture}
\end{center}
\caption{$J_{n}$ as the result of $2^i$ infections on $R_i$}\label{fig:alternativeviews}
\end{figure}
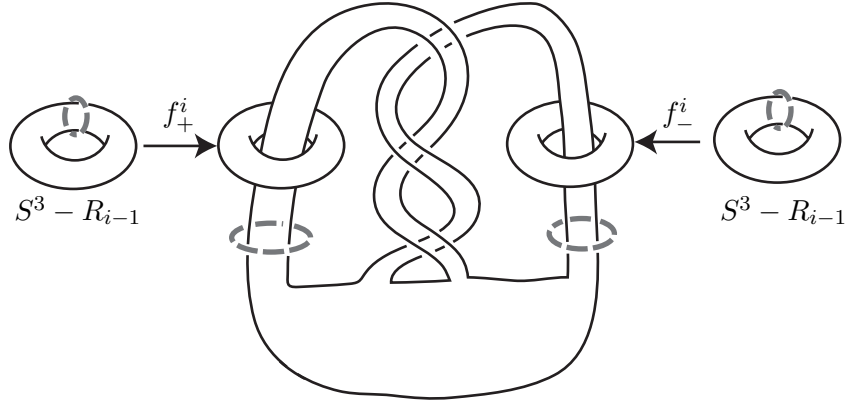
If we alter our point of view by \emph{postponing} (ignoring for the moment) the infections, then we are precisely in the situation of Figure~\ref{fig:mapsf}, that is if we first replace the two fat solid tori by two copies of $S^3-R_{i-1}$ (by convention the maps are named $f_{\pm}^{i}: S^3-R_{i-1}\to S^3-R_{i}$), then we arrive, by definition, at $R_i$. The two collections of images in $S^3-R_i$ of the $2^{i-1}$ clones are precisely the $2^i$ clones $\{\eta^{i}_*\}\equiv ~\{f^{i}_{\pm}\circ\dots\circ f^1_{\pm}(\eta^0)\}$. If we \emph{then} perform these $2^i$ infections using the knot $J_{n-i}(K)$ as the infecting knot in each case, we arrive at the description claimed in the proposition. This completes the inductive step.
\end{proof}

\begin{cor}\label{cor:infection} $J_n(K)$ may be obtained from the ribbon knot $R_n$ as the result of $2^n$ infections along clones, $\{f^{n}_{\pm}\circ\dots\circ f^1_{\pm}(\eta^0)\}$, that lie in $\pi_1(S^3\backslash R_n)^{(n)}$, using the knot $K$ as the infecting knot each time.
\end{cor}

\begin{proof}[Proof of Corollary~\ref{cor:infection}] Apply Proposition~\ref{prop:altdescriptions} in the case $i=n$.
\end{proof}

Returning to the proof of Lemma~\ref{cor:linknsolvable}, suppose that we view the trivial link, $T$, the positive integer $k$ and the curve $\alpha\in F^{(k)}-F^{(k+1)}$ as fixed. Then $T(\alpha,~-)$ may be thought of as an operator from knots to $m$-component links. From this viewpoint, the proof of Proposition~\ref{prop:linkaltdescriptions} below is merely to apply this operator to the result of Proposition~\ref{prop:altdescriptions} above. More details are given below.

\begin{prop}\label{prop:linkaltdescriptions}  For any knot $K$, and any $j,n$ such that $k\leq j\leq n$,  $L_n(K)$ can be obtained from $L_j(U)$ by multiple infections along the $2^{j-k}$ clones $\alpha^{j-k}_*=\{g^{j-k}(\eta_*^{j-k})\}$, using the knot $J_{n-j}(K)$ as the infecting knot in each case, and the clones lie in $\pi_1(S^3-L_j(U))^{(j)}$.
\end{prop}

Assuming Proposition~\ref{prop:linkaltdescriptions} momentarily, Lemma~\ref{cor:linknsolvable} follows immediately. Merely apply Proposition~\ref{prop:linkaltdescriptions} with $j=n$. We claim that $L_n(U)$ is a slice link since it is obtained from the slice link $T$ by infecting using the slice knot $R_{n-k}$ (this is an easy exercise for the reader).
\end{proof}

\begin{proof}[Proof of Proposition~\ref{prop:linkaltdescriptions}] By definition,
$$
L_n(K)\equiv T(\alpha,J_{n-k}(K)), ~L_j(U)\equiv T(\alpha,J_{j-k}(U)).
$$
Since $0\leq j-k\leq n-k$, we have from Proposition~\ref{prop:altdescriptions} that $J_{n-k}(K)$ can be obtained from $J_{j-k}(U)\cong R_{j-k}$ by multiple infections along the $2^{j-k}$ clones $\{\eta^{j-k}_*\}$, using the knot $J_{n-j}(K)$ as the infecting knot in each case. Moreover each clone $\eta^{j-k}_*$ lies in $\pi_1(S^3-R_{j-k})^{(j-k)}$. Therefore, postponing the infections as in Proposition~\ref{prop:altdescriptions}, and as suggested by Figure~\ref{fig:linkaltdescriptions}, we see that $L_n(K)\equiv T(\alpha,J_{n-k}(K))$ can be obtained from $L_j(U)\equiv T(\alpha,R_{j-k})$ by multiple infections along the circles $\{\alpha^{j-k}_*\}=\{g^{j-k}(\eta_*^{j-k})\}$ (that we shall also call \textbf{clones}) using the knot $J_{n-j}(K)$ as the infecting knot in each case.
\begin{figure}[htbp]
\setlength{\unitlength}{1pt}
\begin{center}
\begin{picture}(230,127)
\put(0,0){\includegraphics{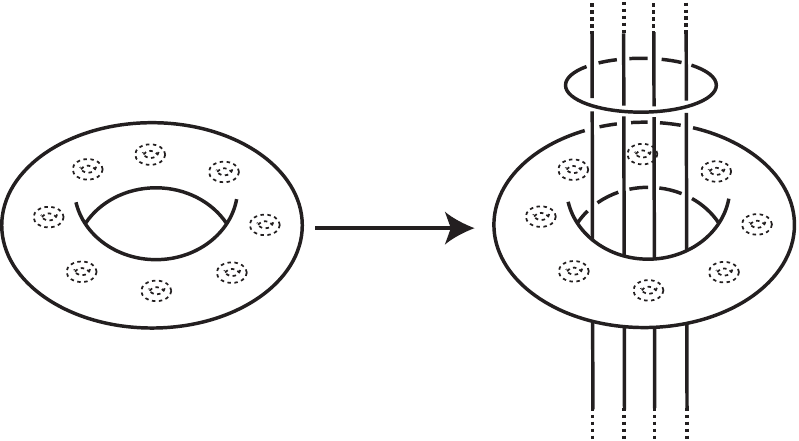}}
\put(148,100){$\alpha^+$}
\put(108,68){$g^{j-k}$}
\put(-22,62){$$}
\put(-22,52){}
\put(-55,62){$\{\eta^{j-k}_*\}\longrightarrow$}
\put(25,19){$S^3-R_{j-k}$}
\put(210,19){$T$}
\end{picture}
\end{center}
\caption{$T(\alpha,J_{n-k}(K)$) obtained from $T(\alpha, R_{j-k})$}\label{fig:linkaltdescriptions}
\end{figure}
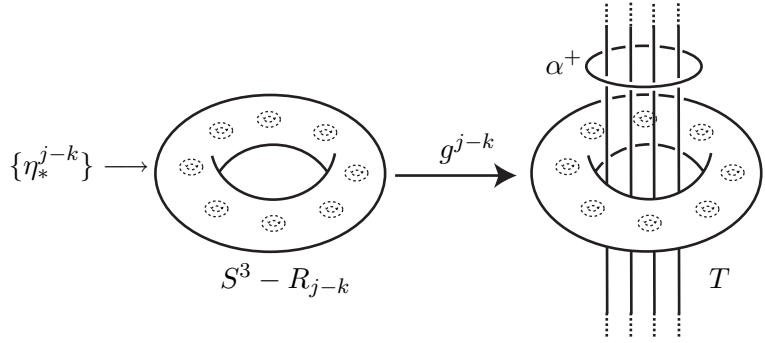

Since $\alpha\in \pi_1(S^3-T)^{(k)}$, the technical result ~\cite[proof of Theorem 8.1]{C} shows that the longitudinal push-off, $\alpha^+$, of $\alpha$ lies in $\pi_1(S^3-\alpha)^{(j)}$ and thus in $\pi_1(S^3-T(\alpha,R_{j-k}))^{(k)}$. Hence, since the meridian of $R_{j-k}$ is identified with $\alpha^+$,
$$
g^{j-k}_*(\pi_1(S^3-R_{j-k}))\subset \pi_1(S^3-L_j(U))^{(k)},
$$
(recalling that $L_j(U)\equiv T(\alpha,J_{n-k}(U))\equiv T(\alpha,R_{n-k})$). Since, by Proposition~\ref{prop:altdescriptions}, each clone $\eta^{j-k}_*$ lies in $\pi_1(S^3-R_{j-k})^{(j-k)}$, each clone $g^{j-k}(\eta^{j-k}_*)$ lies in $\pi_1(S^3-L_j(U))^{(j)}$.

This completes the proof of Proposition~\ref{prop:linkaltdescriptions}.
\end{proof}

Finally we give the proof of Theorem~\ref{thm:iteratedlink}.

\begin{proof}[Proof of Theorem~\ref{thm:iteratedlink}] First we need some notation.

\begin{defn}\label{defn:linkghosts} Let $\mu_j$ denote a meridian of $R_j$ for $0\leq j\leq n-k$. A \textbf{ghost} of $\mu_j$, denoted $(\mu_j)_*$ is an element of the set of $2^{n-k-j}$ circles $\{g^{n-k}f^{n-k}_{\pm}\circ\dots\circ f^{j+1}_{\pm}(\mu_j)\}$. Thus, for any $j$, the ghosts of $\mu_j$ live in $S^3-T(\alpha,R_{n-k})$ and represent elements of $\pi_1(S^3-T(\alpha,R_{n-k}))^{(n-j)}$. These circles are precisely the meridians of the \emph{copies} of $S^3-R_j$ that are embedded in $S^3-T(\alpha,R_{n-k})$ via the maps $\{g^{n-k}f^{n-k}_{\pm}\circ\dots\circ f^{j+1}_{\pm}\}$. Note that $\mu_0$ is the meridian of $R_0=U$ so $\mu_0=\eta^0$. Thus in particular, taking $j=0$, the ghosts of $\mu_0$ coincide with the clones $\{\alpha^{n-k}_*\}$, that is $\{(\mu_0)_*\}=\{\alpha^{n-k}_*\}$.
\end{defn}

Observe that Theorem~\ref{thm:iteratedlink} is a special case ($j=0$) of the following more general result. This proposition should be viewed as a formulation of the inductive step in an inductive proof of Theorem~\ref{thm:iteratedlink}. Hence we may consider the proof of Theorem~\ref{thm:iteratedlink} is finished, but we owe the reader a proof of the following.
\end{proof}

\begin{prop}\label{prop:iteratedghost} Suppose $0\leq j \leq n-k$ and $W$ is an \emph{arbitrary} rational $(n-j)$-solution for $T_{n-k}\equiv T(\alpha,R_{n-k})$. Then at least one of the ghosts of $\mu_{j}$ maps non-trivially under the inclusion-induced map
$$
j_*:\pi_1(M_{T_{n-k}})\to \pi_1(W)/\pi_1(W)_r^{(n-j+1)}.
$$
\end{prop}

\begin{proof}[Proof of Proposition~\ref{prop:iteratedghost}] Here we view $k$ and $n$ as fixed and proceed by downward induction on $j$. First suppose $j=n-k$. In this degenerate case the single ghost is merely the meridian of $R_{n-k}$ viewed as a circle in $T(\alpha, R_{n-k})$, which is of course identified with a push-off, $\alpha^+$, of $\alpha$ itself, and $W$ is a rational $(k)$-solution for $M_{T_{n-k}}$. We must show that $j_*(\alpha^+)\neq 1$ under the map
$$
j_*:\pi_1(M_{T_{n-k}})\to \pi_1(W)/\pi_1(W)^{(k+1)}_r.
$$

Since $T_{n-k}$ is obtained from the trivial link $T$ by infection on a curve $\alpha\in F^{(k)}$, by ~\cite[Proposition 3.1]{Lei3}, there is a degree one map $r:M_{T_{n-k}}\to M_T$ that induces an isomorphism
$$
\pi_1(M_{T_{n-k}})/(\pi_1(M_{T_{n-k}}))^{(k+1)}_r\cong F/F^{(k+1)}
$$
and sends $\alpha^+$ to $\alpha$. Since $\alpha$ is not in  $F^{(k+1)}$, $\alpha^+ \neq 1$ in $\pi_1(M_{T_{n-k}})/\pi_1(M_{T_{n-k}})^{(k+1)}_r$. This also implies that the successive terms of the derived series of $\pi_1(M_{T_{n-k}})$ agree with those of the free group (up to this value of $k$). Thus the derived series, the rational derived series and even Harvey's torsion-free derived series agree for this group (up to this value of $k$)~\cite[Section 2]{Ha2})~\cite[Proposition 2.3]{Ha2}. This is useful because we now claim that the following is a monomorphism
$$
\pi_1(M_{T_{n-k}})/\pi_1(M_{T_{n-k}})^{(k+1)}_r\overset{j_*}\to \pi_1(W)/\pi_1(W)^{(k+1)}_r
$$
because the composition
$$
\pi_1(M_{T_{n-k}})/\pi_1(M_{T_{n-k}})^{(k+1)}_r\overset{j_*}\to \pi_1(W)/\pi_1(W)^{(k+1)}_r\to \pi_1(W)/\pi_1(W)^{(k+1)}_H
$$
is a monomorphism by the following result of the authors. Here we use that $W$ is a rational $(k)$-solution for $M_{T_{n-k}}$ and that the torsion-free derived series of a free group is the same its rational derived series.

\begin{prop}[Proposition 4.11 ~\cite{CH2}]\label{prop:nsolvable} If $M$ is rationally $(k)-solvable$ via $W$ then, letting $A=\pi_{1}(M)$ and $B=\pi_{1}(W)$, the inclusion $j:M \rightarrow W$ induces a monomorphism
$$j_*: \frac{\pi_1(M)}{\pi_1(M)_{H}^{(k+1)}} \hookrightarrow  \frac{\pi_1(W)}{\pi_1(W)_{H}^{(k+1)}}.$$
\end{prop}

\noindent It follows that $j_*(\alpha^+)\neq 1$ as required by Proposition~\ref{prop:iteratedghost}. Thus the Proposition holds for $j=n-k$.

Now suppose that the Proposition is true for $j+1$ where $1\leq j+1 \leq n-k$. We will establish it for $j$ (downwards induction). So consider a rational $(n-j)$-solution, $W$, for $M_{T_{n-k}}$. Let $\Lambda=\pi_1(W)/\pi_1(W)^{(n-j)}_r$ and let $\psi:\pi_1(W)\to \Lambda$, and $\phi:\pi_1(M_{T_{n-k}})\to \Lambda$ be the induced coefficient systems. Note that $W$ is \emph{a fortiori} a rational ($n-j-1$)-solution. Therefore the inductive hypothesis applies to $W$ for the value $j+1$ and allows us to conclude that at least one ghost of $\mu_{j+1}$ does not map into $\pi_1(W)^{(n-j)}_r$ under the inclusion, that is, we have ~$\phi((\mu_{j+1})_*)\neq 1$ for some ghost of $\mu_{j+1}$. We will need this fact below.

We can apply Proposition~\ref{prop:linkaltdescriptions} with $K=U$  to deduce that $L_n(U)$ ($\equiv T(\alpha,R_{n-k})\equiv T_{n-k}$) can be obtained from $L_{n-j-1}(U)\equiv T_{n-j-k-1}$ by infections along the clones $\{\alpha^{n-j-k-1}_*\}= \{g^{n-j-k-1}(\eta_*^{n-j-k-1})\}$  using the knot $R_{j+1}$ as infecting knot in each case. Then, in the notation of Theorem~\ref{thm:nontriviality}
$$
T_{n-k}=T_{n-j-k-1}(\alpha^{n-k-j-1}_i,R_{j+1}^i, 1\leq i\leq 2^{n-k-j-1})
$$
where $(R_{j+1})^i$ is the $i^{th}$ copy of $R_{j+1}$. Applying Theorem~\ref{thm:nontriviality} we see that, for any clone such that $\phi((\alpha^{n-k-j-1}_i)^+)\neq 1$ the kernel, $P_i$ of the composition
$$
\mathcal{A}_0(R_{j+1})\overset{}\to (\mathcal{A}_0(R_{j+1}) \otimes\mathbb{Q}\Lambda)\overset{i_*}{\to} H_1(M_{T_{n-k}};\mathbb{Q}\Lambda)\overset{j_*}\to H_1(W;\mathbb{Q}\Lambda),
$$
satisfies $P_i\subset P_i^\perp$. We claim that there exists at least one such clone. For, by definition of infection, when we infect $T_{n-j-k-1}$ along $\alpha^{n-k-j-1}_i$ the push-off or longitude of such a circle, $(\alpha^{n-k-j-1}_i)^+$, is identified to the meridian of the $i^{th}$ copy of the infecting knot $(R_{j+1})^i$. This meridian, when viewed as a circle in $T_{n-k}$, is not a meridian of the abstract knot $R_{j+1}$, but rather an embedded copy of that meridian in $T_{n-k}$. Thus $(\alpha^{n-k-j-1}_i)^+$, viewed as a circle in $T_{n-k}$, is, by definition, one of the \emph{ghosts} of $\mu_{j+1}$! But we established above, by our inductive assumption, that for at least one of these ghosts, $\phi((\mu_{j+1})_*)\neq 1$. Thus we have verified that there is at least one such clone (say the $i^{th}$) for which the hypotheses of Theorem~\ref{thm:nontriviality} apply. We now restrict attention to such a value of $i$.

The two circles
$$
f^{j+1}_\pm(\mu_{j}) \in \pi_1(S^3-R_{j+1})^{(1)}
$$
as shown in the Figure~\ref{fig:twocircles}, form a generating set for $\mathcal{A}_0(R_{j+1})$ (which is isomorphic to $\mathcal{A}_0(R_{1})$ and hence nontrivial).
\begin{figure}[htbp]
\setlength{\unitlength}{1pt}
\begin{center}
\begin{picture}(377,181)
\put(0,0){\includegraphics{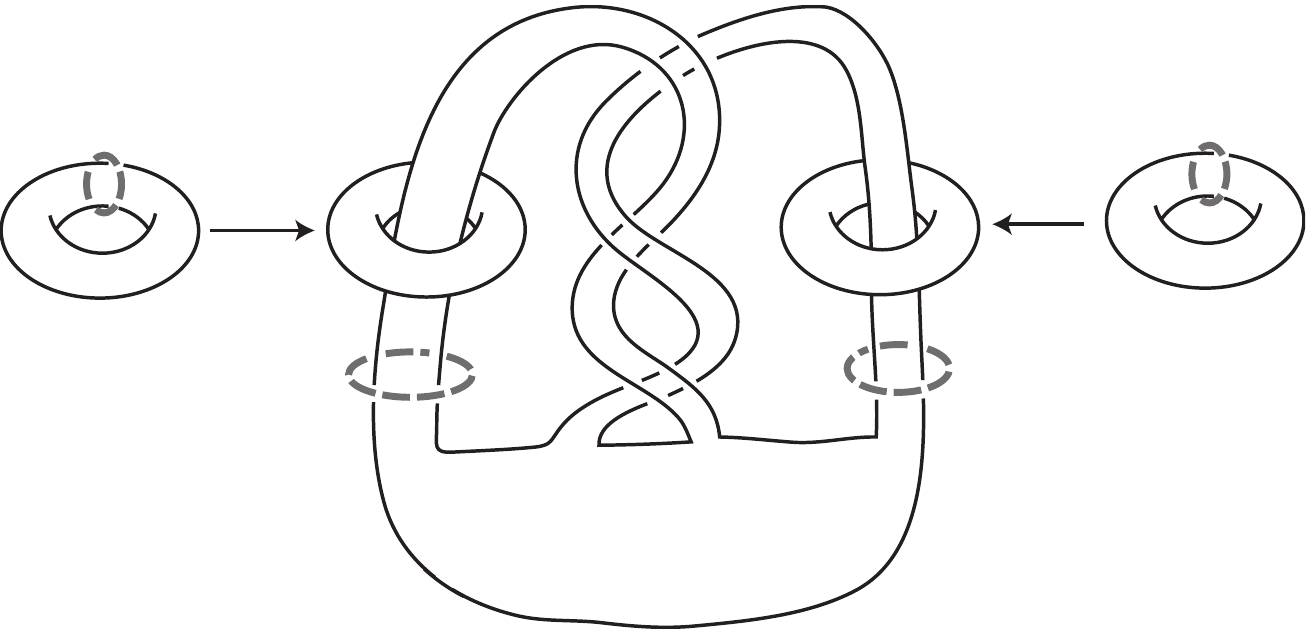}}
\put(67,123){$f^{j+1}_+$}
\put(293,123){$f^{j+1}_-$}
\put(10,84){$S^3-R_{j}$}
\put(330,86){$S^3-R_{j}$}
\put(63,64){$f^{j+1}_+(\mu_{j})$}
\put(275,64){$f^{j+1}_-(\mu_{j})$}
\put(19,142){$\mu_{j}$}
\put(334,143){$\mu_{j}$}
\end{picture}
\end{center}
\caption{Inside the $i^{th}$ copy of $S^3-R_{j+1}$}\label{fig:twocircles}
\end{figure}
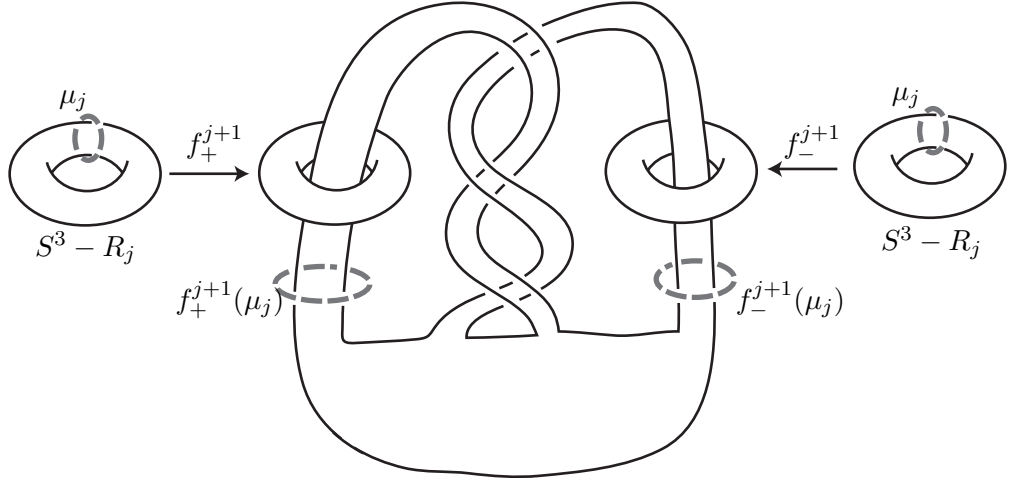

From this we can conclude that at least one of the generators is not in $P_i$ since otherwise
$$
P_i=\mathcal{A}_0(R_{j+1})\subset\mathcal{A}_0(R_{j+1})^\perp,
$$
contradicting the nonsingularity of the classical Blanchfield form of $\mathcal{A}_0(R_{j+1})$. Finally, consider the commutative diagram below, where we abbreviate $\pi_1(W)$ by $\pi$. Recall that $H_1(W;\mathbb{Q}\Lambda)$ is identifiable as the ordinary rational homology of the covering space of $W$ whose fundamental group is the kernel of $\psi:\pi\to \Lambda$. Since this kernel is precisely $\pi^{(n-j)}_r$, we have that
$$
H_1(W;\mathbb{Q}\Lambda)\cong (\pi^{(n-j)}_r/[\pi^{(n-j)}_r,\pi^{(n-j)}_r])\otimes_{\mathbb{Z}} \mathbb{Q}
$$
as indicated in the diagram below. By the definition of the rational derived series, the far-right vertical map $j$ is injective.
\begin{equation*}
\begin{CD}
\pi_1(S^3-R_{j+1})^{(1)}      @>i_*>>    \pi_1(M_{T_{n-k}})^{(n-j)}  @>j_*>>   \pi^{(n-j)}_r  @>>>
\frac{\pi^{(n-j)}_r}{\pi^{(n-j+1)}_r} \\
  @VVV   @VVV        @VVV       @VVjV\\
\mathcal{A}_0(R_{j+1})     @>i_*>>  H_1(M_{T_{n-k}};\mathbb{Q}\Lambda)    @>j_*>>  H_1(W;\mathbb{Q}\Lambda) @>\cong>>
  \frac{\pi^{(n-j)}_r}{[\pi^{(n-j)}_r,\pi^{(n-j)}_r]}\otimes_{\mathbb{Z}} \mathbb{Q}\\
\end{CD}
\end{equation*}
Hence, since the composition in the bottom row sends one of the two homology classes $[f^{j+1}_\pm(\mu_{j})]$ to non-zero, the composition in the top row sends at least one of the two $f^{j+1}_\pm(\mu_{j})$ to non-zero under $i_*$. Now observe that the map $i_*$ in the top row above is induced by one of the compositions $g^{n-k}\circ f_\pm^{n-k}\circ\dots\circ f_\pm^{j+2}$. Thus
$$
i_*(f^{j+1}_\pm(\mu_{j}))=g^{n-k}\circ f_\pm^{n-k}\circ\dots\circ f_\pm^{j+2}\circ f^{j+1}_\pm(\mu_{j}).
$$
For various values of $\pm$ these are precisely the ghosts of $\mu_{j}$. Hence we have shown that for at least one such ghost of $\mu_{j}$
$$
j_*((\mu_{j})_*)\neq 1 ~\text{in} ~\pi^{(n-j)}_r/\pi^{(n-j+1)}_r
$$
as desired.

This finishes the inductive proof of Proposition~\ref{prop:iteratedghost}.
\end{proof}

Since we did not use very heavily the fact that $T$ is a trivial link nor did we use much about the specific nature of the ribbon knot $9_{46}$,  the proof shows the following more general result.

\begin{thm}\label{thm:mainlink3} Suppose $T$ is a slice link, $\alpha$ is an unknotted circle in $S^3-T$ that represents an element in $\pi_1(S^3-T)^{(k)}$ but not in $\pi_1(M_T)^{(k+1)}_H$. Suppose for each $j$, $1\leq j\leq n-k$, $R_j$ is a slice knot, $\{\eta_{j1},\dots,\eta_{jm_j}\}$ is a trivial link of circles in $S^3-R_j$ with the property that the submodule of the classical Alexander polynomial of $R_j$ generated by $\{\eta_{j1},\dots,\eta_{jm_j}\}$ contains elements $x,y$ such that $\mathcal{B}\ell_0^j(x,y)\neq 0$, where $\mathcal{B}\ell_0^j$ is the Blanchfield form of $R_j$. Finally suppose that Arf($K$)$=0$. Then the result, $L(K)\equiv T_{\alpha}\circ R_{n-k}\circ\dots\circ R_1(K)$, of the iterated generalized doubling (applied to $K$) lies in $\mathcal{F}_{n}$ and there is a constant $C$ (independent of $K$), such that if $|\rho_0(K)|>C$, then $L(K)$ is of infinite order in the topological concordance group (moreover no multiple lies in $\mathcal{F}_{n+1}$).
\end{thm}

\section{Higher-Order Signatures as Obstructions to being Slice and the COT $(n)$-solvable Filtration }\label{sec:appendix}

%\textbf{The COT $\boldsymbol{n}$-solvable filtration}

Recall that \cite[Section 8]{COT} introduced a filtration of the concordance classes of links $\mathcal{C}$
$$
\cdots \subseteq \mathcal{F}_{n} \subseteq \cdots \subseteq
\mathcal{F}_1\subseteq \mathcal{F}_{0.5} \subseteq \mathcal{F}_{0} \subseteq \mathcal{C}.
$$
where the elements of $\mathcal{F}_{n}$ and $\mathcal{F}_{n.5}$ are called \emph{$(n)$-solvable links} and \emph{$(n.5)$-solvable links} respectively. In the case of knots this is a filtration by \emph{subgroups} of the knot concordance group. A slice link $L$ has the property that its zero surgery $M_L$ bounds a $4$-manifold $W$ (namely the exterior of the slicing disks) such that $H_1(M_L)\to H_1(W)$ is an isomorphism and $H_2(W)=0$. An \emph{$(n)$-solvable} link is one, loosely speaking, such that $M_L$ bounds a $4$-manifold $W$ such that $H_1(M_L)\to H_1(W)$ is an isomorphism and the intersection form on $H_2(W)$ ``looks'' hyperbolic modulo the $n^{th}$-term of the derived series of $\pi_1(W)$. We shall only give a detailed definition of the slightly larger class of \emph{rationally $(n)$-solvable links}.

For a compact oriented topological 4-manifold $W$, let $W^{(n)}$ denote the covering
space of $W$ corresponding to the $n$-th derived subgroup of $\pi_1(W)$. The deck translation group of this cover is the solvable group $\pi_1(W)/\pi_1(W)^{(n)}$. Then $H_2(W^{(n)};\mathbb{Q})$ can be endowed with the structure of a right $\mathbb{Q}[\pi_1W)^{(n)}]$-module. This agrees with the homology group with twisted coefficients $H_2(W;\mathbb{Q}[\pi_1(W)^{(n)}])$. There is an equivariant
intersection form
$$
\lambda_n : H_2(W^{(n)};\mathbb{Q}) \times H_2(W^{(n)};\mathbb{Q}) \lra
\mathbb{Q}[\pi_1(W)/\pi_1(W)^{(n)}]
$$
 \cite[Section 7]{COT}\cite[Chapter 5]{Wa}. The usual intersection form is the case $n=0$. In general, these
intersection forms are singular. Let $I_n \equiv$ image($j* : H_2(\partial W^{(n)};\mathbb{Q}) \to H_2(W^{(n)};\mathbb{Q})$). Then this intersection form factors
through
$$
\ov{\lambda_n} : H_2(W^{(n)};\mathbb{Q})/I_n \times H_2(W^{(n)};\mathbb{Q})/I_n \lra \mathbb{Q}[\pi_1(W)/\pi_1(W)^{(n)}].
$$
We define a  \textbf{rational $\boldsymbol{(n)}$-Lagrangian} of $W$ to be a
submodule of $H_2(W;\mathbb{Q}[\pi_1W)^{(n)}]$ on which
$\ov\lambda_n$ vanishes identically and which maps onto a $\frac12$-rank
subspace of $H_2(W;\mathbb{Q})/I_0$ under the covering map. An
\textbf{$\boldsymbol{(n)}$-surface} is a based and immersed surface
in $W$ that can be lifted to $W^{(n)}$. Observe that any class in
$H_2(W^{(n)})$ can be represented by an $(n)$-surface and that
$\lambda_n$ can be calculated by counting intersection points in
$W$ among representative $(n)$-surfaces weighted appropriately by
signs and by elements of $\pi_1(W)/\pi_1(W)^{(n)}$. We say a rational
$(n)$-Lagrangian $L$ admits \textbf{rational $\boldsymbol{(m)}$-duals} (for $m\le n$) if $L$
is generated by (lifts of) $(n)$-surfaces $\ell_1,\ell_2,\ldots,\ell_g$ and
there exist $(m)$-surfaces $d_1,d_2,\ldots, d_g$ such that $H_2(W;\mathbb{Q})/I_0$
has rank $2g$ and $\lambda_m(\ell_i,d_j)=\delta_{i,j}$.

Under the assumption that we will impose below, that
$$
H_1(M;\mathbb{Q})\to H_1(W;\mathbb{Q})
$$
is an isomorphism, it follows that the dual map
$$
H_3(W,M;\mathbb{Q})\to H_2(M;\mathbb{Q})
$$
is an isomorphism and hence that $I_0=0$. Thus the ``size'' of rational $(n)$-solutions is dictated by the rank of $H_2(W;\mathbb{Q})$.

\begin{defn}
\label{defn:rationalnsolvable} Let $n$ be a nonnegative integer. A compact, connected oriented topological 4-manifold $W$ with $\partial W = M$ is a \textbf{rational $\boldsymbol{(n)}$-solution for $\boldsymbol{M}$} if
\begin{itemize}
\item $H_1(M;\mathbb{Q})\to H_1(W;\mathbb{Q})$ is an isomorphism, and
\item  $W$ admits a rational $(n)$-Lagrangian with rational $(n)$-duals.
\end{itemize}
Then we say that \textbf{$\boldsymbol{M}$ is rationally $\boldsymbol{(n)}$-solvable via $\boldsymbol{W}$}. A link $L$ is a \textbf{rationally $\boldsymbol{(n)}$-solvable link} if
$M_L$ is rationally $(n)$-solvable for some such $W$.
\end{defn}

\begin{defn}
\label{defn:rationaln.5solvable} Let $n$ be a nonnegative integer. A compact, connected oriented 4-manifold $W$ with
$\partial W = M$ is a \textbf{rational $\boldsymbol{(n.5)}$-solution for $\boldsymbol{M}$} if
\begin{itemize}
\item $H_1(M;\mathbb{Q})\to H_1(W;\mathbb{Q})$ is an isomorphism, and
\item  $W$ admits a rational $(n)$-Lagrangian with rational $(n+1)$-duals.
\end{itemize}
Then we say that \textbf{$\boldsymbol{M}$ is rationally $\boldsymbol{(n.5)}$-solvable via $\boldsymbol{W}$}. A link $L$ is a \textbf{rationally $\boldsymbol{(n.5)}$-solvable link} if $M_L$ is rationally $(n.5)$-solvable for some such $W$.
\end{defn}

A $4$-manifold $W$ is an \textbf{$\boldsymbol{(n)}$-solution} (respectively an \textbf{$\boldsymbol{(n.5)}$-solution}) if, in addition, it is spin, it satisfies the conditions above with $\mathbb{Q}$ replaced by $\mathbb{Z}$ and the equivariant self-intersection form also vanishes on the Lagrangian (see ~\cite[Section 8]{COT}.

\begin{remark}
\

\label{rem:n-solvable}
\begin{enumerate}
\item An $(n)$-solution is a fortiori a rational $(n)$-solution.
\item An $(n)$-solution (respectively rational $(n)$-solution) is a fortiori an $(m)$-solution (respectively rational $(m)$-solution) for any $m<n$.
\item If $L$ is slice in a topological (rational) homology $4$-ball then the complement of a set of slice disks is an (rational) $(n)$-solution for any integer or half-integer $n$.  This follows since if $H_2(W;\mathbb{Z})=0$ then the Lagrangian may be taken to be the zero submodule.
\end{enumerate}
\end{remark}

The following result is useful.

\begin{lem}\label{lem:nsolv} Suppose $L$ is a link obtained from a $(p+q)$-solvable link $R$ by infection along curves in $\pi_1(S^3-R)^{(p)}$ using knots $K_i$. Suppose the knots $K_i$ are $(q)$-solvable. Then $L$ is also a $(p+q)$-solvable link.
\end{lem}
\begin{proof} One can repeat almost verbatim the proof of ~\cite[Proposition 3.1]{COT2} (see also ~\cite[Corollary 3.14]{CT}). However, one also needs the following result.
\begin{lem}\label{lem:normgender} Suppose $\phi:A\to B$ is a group homomorphism that is surjective on abelianizations. Then, for any positive integer $n$, $\phi(A)$ normally generates $B/B^{(n)}$.
\end{lem}
\begin{proof}[Proof of Lemma~\ref{lem:normgender}] The proof is by induction on $n$. The case $n=1$ is the hypothesis. Now consider $b\in B$. Then $b=\phi(a)\prod_{i=1}^m[b_{i1},b_{i2}]$ where $a\in A$ and $b_{i1},b_{i2}\in B$.  It now suffices to show that a single commutator
$$
[b_1,b_2]\in  <\phi(A)>B^{(n)}
$$
where $<\phi(A)>$ denotes the normal closure in $B$. By the inductive hypothesis
$$
b_j\in <\phi(A)>B^{(n-1)}
$$
for $j=1,2$. Hence
$$
[b_1,b_2]\in  [<\phi(A)>B^{(n-1)},<\phi(A)>B^{(n-1)}],
$$
which equals $<\phi(A)>B^{(n)}$ by simple commutator calculus.
\end{proof}
This completes the proof of Lemma~\ref{lem:nsolv}.
\end{proof}

\begin{thm}[Cochran-Orr-Teichner~{\cite[Theorem 4.2]{COT}}]\label{thm:sliceobstr} If a knot $K$ is rationally $(n.5)$-solvable via $W$ and  $\phi:\pi_1(M_K)\to \G$ is a PTFA coefficient system that extends to $\pi_1(W)$ and such that $\G^{(n+1)}=1$, then $\rho(M_K,\phi)=0$.
\end{thm}

For links the following recent result of the first two authors is the best known result (although see ~\cite[Theorem 5.9]{CHL3}). Note the extra rank condition.

\begin{thm}[Cochran-Harvey~{\cite[Theorem 4.9, Proposition 4.11]{CH2}}]\label{thm:rho=0} Let $\Gamma$ be a PTFA group such that $\G^{(n+1)}=0$. Let $M$ be a closed, connected, oriented $3$-manifold equipped with a non-trivial coefficient system $\phi:\pi_1(M)\to \Gamma$. Suppose $\text{rank}_{\mathcal{K}\Gamma}(H_1(M;\mathcal{K}\G))= \beta_1(M)-1$. Then if $M$ is
rationally
$(n.5)$-solvable via a
$4$-manifold $W$ over which $\phi$ extends, then
$$
\rho(M,\phi)= \sigma^{(2)}_\Gamma(W)-\sigma(W)=0.
$$
Moreover, if additionally $M$ is
rationally $(n+1)$-solvable via $W$ then the extra rank condition above is automatically satisfied.
\end{thm}

\begin{proof}[Proof that Theorem~\ref{thm:rho=0} implies Theorem~\ref{thm:linksliceobstr}] Since $\G$ is PTFA, it is solvable so there exists some $n$ such that $\G^{(n+1)}=0$. Let $W$ denote the exterior of the slicing disks. By Alexander duality, $H_2(W;\mathbb{Q})=0$ and $H_1(M_L;\mathbb{Q})\to H_1(W;\mathbb{Q})$ is an isomorphism. Thus $W$ is a certainly a rational $(n+1)$-solution for $L$. Then the result follows immediately from Theorem~\ref{thm:rho=0}.
\end{proof}

There is another common situation in which the extra rank condition is satisfied.

\begin{lem}\label{lem:rank} Suppose $L$ is a link obtained from the link $R$ by infections on circles $\eta_i$ using knots $K_i$. Suppose $\phi:\pi_1(M_L)\to\G$ is a nontrivial PTFA coefficient system such that $\phi(\mu_{\eta_i}\equiv l_{K_i})=1$. Then there is a coefficient system $\phi:\pi_1(M_L)\to\G$ induced on $M_R$ and
$$
\text{rank}_{\mathcal{K}\G}(H_1(M_{L};\mathcal{K}\G))\geq \text{rank}_{\mathcal{K}\G}(H_1(M_R;\mathcal{K}\G)).
$$
In particular if $R$ is the trivial link of $m$ components then
$$
\text{rank}_{\mathcal{K}\G}(H_1(M_{L};\mathcal{K}\G))= \beta_1(M_L)-1.
$$
\end{lem}
\begin{proof}[Proof of Lemma~\ref{lem:rank}] Consider the cobordism $E_L$ of Figure~\ref{fig:mickey}. By Property $(1)$ of Lemma~\ref{lem:mickeyfacts}, the map
$$
\pi_1(M_L)\to \pi_1(E_L)
$$
is a surjection whose kernel is normally generated by $\{\mu_{\eta_i}\}$. Thus, as shown there, $\phi$ extends uniquely to $\pi_1(E_L)$ and hence by restriction to $\pi_1(M_R)$. Therefore there is a surjection
$$
H_1(M_L;\mathcal{K}\G)\to H_1(E_L;\mathcal{K}\G)
$$
so
$$
\text{rank}_{\mathcal{K}\G}(H_1(M_L;\mathcal{K}\G))\geq \text{rank}_{\mathcal{K}\G}(H_1(E_L;\mathcal{K}\G)).
$$
Now examine the Mayer-Vietoris sequence with $\mathcal{K}\G$ coefficients for $E_L$ as in the proof of Lemma~\ref{lem:mickeysig}
$$
\oplus_i H_1(\eta_i \times D^2)\to \oplus_i H_1(M_{K_i})\oplus H_1(M_R)\to H_1(E_L)\overset{\partial_*}\to \oplus_iH_0(\eta_i \times D^2).
$$
We claim that the inclusion-induced maps
$$
H_0(\eta_i\x D^2;\mathcal{K}\G)\ra H_0(M_i;\mathcal{K}\G)
$$
are injective. In the case that $\phi(\eta_i)\neq 1$, $H_0(\eta_i\x D^2;\mathcal{K}\G)=~0$ by ~\cite[Proposition 2.9]{COT}, so injectivity holds. If $\phi(\eta_i)=1$ then, since $\eta_i$ is equated to the meridian of $K_i$, $\phi(\mu_{K_i})=1$. Since $\mu_i$ normally generates $\pi_1(M_i)$, it follows that the coefficient systems on $\eta_i\x D^2$ and $M_i$ are trivial and hence the injectivity follows from the injectivity with $\mathbb{Z}$-coefficients, which is obvious since both are path-connected. Hence $\partial_*$ is the zero map.
Similarly we claim that the inclusion-induced maps
$$
H_1(\eta_i\x D^2;\mathcal{K}\G)\ra H_1(M_{K_i};\mathcal{K}\G)
$$
are isomorphisms. In the case that $\phi(\eta_i)\neq 1$, both groups are zero by ~\cite[Lemma 2.10]{COT}. If $\phi(\eta_i)=1$ then both coefficient systems are trivial and result follows from the result for $\mathbb{Z}$-coefficients, which is obvious since $u_{K_i}$ generates $H_1(M_{K_i})\cong \mathbb{Z}$.

Armed with these observations, it now follows from the Mayer-Vietoris sequence that
$$
H_1(M_{R};\mathcal{K}\G)\cong H_1(E_L;\mathcal{K}\G).
$$
and the first result follows.

If $R$ is a trivial link then $\pi_1(M_R)$ is the free group $F$ of rank $m$. But it is easy to see from an Euler characteristic argument \cite[Lemma 2.12]{COT}) that
$$
\text{rank}_{\mathcal{K}\G}(F;\mathcal{K}\G))=\beta_1(F)-1=m-1.
$$
Thus
$$
\text{rank}_{\mathcal{K}\G}(H_1(M_L;\mathcal{K}\G))\geq \beta_1(M_L)-1
$$
but by ~\cite[Proposition 2.11]{COT}, this is also the maximum this rank can achieve, so the inequality is an equality.
\end{proof}

\bibliographystyle{plain}
\bibliography{mybib5mathscinet}
\end{document}